\documentclass[hidelinks,onefignum,onetabnum]{siamart220329}

\usepackage[T1]{fontenc}
\usepackage{lipsum}
\usepackage{amsfonts}
\usepackage{graphicx}
\usepackage{epstopdf}
\usepackage{siunitx}
\usepackage{amsopn}
\usepackage{amsmath}

\DeclareMathOperator*{\argmin}{arg\,min}
\usepackage{algorithmic}
\usepackage{color}
\usepackage{pifont}
\usepackage{caption}
\usepackage{subcaption}
\usepackage{pgfplots}
\pgfplotsset{compat=newest}
\usepackage{tikz}
\usetikzlibrary{arrows.meta,calc,3d,shapes,math,fit,positioning,shadows,patterns,backgrounds}
\usepgfplotslibrary{groupplots,colorbrewer,colormaps}
\tikzset{no shadows/.style={general shadow/.style=}}
\def\myscale{0.5}
\tikzset{block/.style={rectangle split, draw, rectangle split parts=2,
		text width=14em*\myscale, text centered, rounded corners, minimum height=4em*\myscale},
grnblock/.style={materia, fill=green!20, text width=10em*\myscale, text centered, rounded corners, minimum height=4em*\myscale},
materia/.style={draw, fill=white!20, text width=15.0em*\myscale, text centered, minimum height=1.5em*\myscale,drop shadow},
whtblock/.style={no shadows, materia, text width=30em*\myscale, minimum width=20em*\myscale, minimum height=10em*\myscale, rounded corners, drop shadow},
line/.style={draw, -{Latex[length=2mm*\myscale,width=1mm]}},
cloud/.style={draw, ellipse,fill=white!20, node distance=3cm*\myscale,    minimum height=4em*\myscale},
container/.style={draw, rectangle,dashed,inner sep=0.9cm*\myscale, rounded
		corners,fill=orange!10,minimum height=2cm*\myscale},
container2/.style={draw, rectangle,dashed,inner sep=0.9cm*\myscale, rounded
		corners,fill=blue!10,minimum height=2cm*\myscale},
container3/.style={draw, rectangle,dashed,inner sep=0.9cm*\myscale, rounded
		corners,fill=green!10,minimum height=2cm*\myscale},
box/.style={draw, rectangle, rounded corners, thick, node
		distance=7em*\myscale,
		text width=6em*\myscale, text centered, minimum height=3.5em*\myscale},
every node/.style={font=\scriptsize}}
\tikzstyle{arrow} = [thick,->,>=stealth]
\pgfplotsset{
	colormap/Set2-5,
	cycle list/Set2-5,
	/pgfplots/layers/Bowpark/.define layer set={
			axis background,axis grid,axis ticks,axis lines,main,axis tick labels,
			axis descriptions,axis foreground
		}{/pgfplots/layers/standard},
	every axis/.append style={
			set layers=Bowpark,
		},
}

\def\minwidth{1}
\def\minheight{0.7}
\tikzstyle{startstop} = [rectangle, rounded corners, 
minimum width=1cm*\minwidth, 
minimum height=1cm*\minheight,
text centered, 
draw=black, 
fill=red!30]

\tikzstyle{io} = [trapezium, 
trapezium stretches=true, 
trapezium left angle=70, 
trapezium right angle=110, 
minimum width=1cm*\minwidth, 
minimum height=1cm*\minheight, text centered, 
draw=black, fill=blue!30]

\tikzstyle{process} = [rectangle, rounded corners,
minimum width=1cm*\minwidth, 
minimum height=1cm*\minheight, 
text centered, 
draw=black, 
fill=orange!30]

\tikzstyle{decision} = [diamond, 
minimum width=1cm*\minwidth, 
minimum height=1cm*\minheight, 
text centered, 
draw=black, 
fill=green!30]
\tikzstyle{arrow} = [thick,->,>=stealth]

\usepackage{mathtools}
\usepackage{listings}
\usepackage[nolist]{acronym}
\usepackage{hyphenat}
\usepackage{enumitem}
\usepackage{verbatim}

\ifpdf
  \DeclareGraphicsExtensions{.eps,.pdf,.png,.jpg}
\else
  \DeclareGraphicsExtensions{.eps}
\fi

\def\Cpp{{C\nolinebreak[4]\hspace{-.05em}\raisebox{.4ex}{\tiny\bf ++}}}

\newsiamremark{remark}{Remark}
\newsiamremark{hypothesis}{Hypothesis}
\crefname{hypothesis}{Hypothesis}{Hypotheses}
\newsiamthm{claim}{Claim}

\headers{Two-level nonlinear Schwarz methods - a parallel implementation}{Kyrill Ho, Axel Klawonn, and Martin Lanser}

\title{Two-level nonlinear Schwarz methods - a parallel implementation  with application to nonlinear elasticity and incompressible flow problems\thanks{{\bf Funding: }{The authors gratefully acknowledge the financial support by the German Federal Ministry of Research, Technology and Space (BMFTR) in the program SCALEXA.
The scientific support and HPC resources provided by the Erlangen National High Performance Computing Center (NHR@FAU) of the Friedrich-Alexander-Universit\"at Erlangen-N\"urnberg (FAU), Germany under the NHR project k107ce is gratefully acknowledged. NHR funding is provided by federal and Bavarian state authorities. NHR@FAU hardware is partially funded by the German Research Foundation (DFG) - 440719683.
}}}

\author{Kyrill Ho\footnotemark[2] \footnotemark[3]%
\and Axel Klawonn\footnotemark[2] \footnotemark[3]%
\and Martin Lanser\footnotemark[2] \footnotemark[3]}

\ifpdf
\hypersetup{
  pdftitle={A Parallel Two-level Nonlinear Schwarz Implementation Applied to 2D Hyperelasticity and Incompressible Flow},
  pdfauthor={Kyrill Ho and Axel Klawonn and Martin Lanser}
}
\fi
\renewcommand{\ttdefault}{cmtt}  
\DeclareTextFontCommand{\mytexttt}{\ttfamily\hyphenchar\font=45\relax}
\makeatletter
\DeclareRobustCommand\ttfamily
        {\not@math@alphabet\ttfamily\mathtt
         \fontfamily\ttdefault\selectfont\hyphenchar\font=-1\relax}
\makeatother
\DeclareTextFontCommand{\mytexttt}{\ttfamily\hyphenchar\font=45\relax}

\def\apply{\mytexttt{apply()}}
\def\nonlinop{\mytexttt{NonLinearOperator}}
\def\nonlinschwarzop{\mytexttt{NonLinearSchwarzOperator}}
\def\simpleoverlappingop{\mytexttt{SimpleOverlappingOperator}}
\def\coarsenonlinschwarzop{\mytexttt{CoarseNonLinearSchwarzOperator}}
\def\schwarzop{\mytexttt{Sch\-warz\-Oper\-ator}}
\def\simplecoarseop{\mytexttt{SimpleCoarseOperator}}
\def\coarseop{\mytexttt{CoarseOperator}}
\def\ipouop{\mytexttt{IPOUHarmonicCoarseOperator}}
\def\overlappingop{\mytexttt{OverlappingOperator}}
\def\combineop{\mytexttt{CombineOperator}}
\def\nonlincombineop{\mytexttt{NonLinearCombineOperator}}
\def\linsumop{\mytexttt{SumOperator}}
\def\nonlinsumop{\mytexttt{NonLinearSumOperator}}

\def\KH{\color{black}}

\begin{document}

\maketitle

{
  \renewcommand{\thefootnote}{\fnsymbol{footnote}}%
  \footnotetext[2]{Department of Mathematics and Computer Science, University of Cologne,
    Weyertal 86-90, 50931 Cologne, Germany,
    (\url{http://www.numerik.uni-koeln.de}),
    (\email{kyrill.ho@uni-koeln.de}, \email{axel.klawonn@uni-koeln.de},
    \email{martin.lanser@uni-koeln.de}).}
  \footnotetext[3]{Center for Data and Simulation Science,
    University of Cologne, Albertus-Magnus-Platz, 50923 Cologne, Germany, (\url{http://www.cds.uni-koeln.de}).}
}

\begin{acronym}
  \acro{FEDDLib}{Finite Element and Domain Decomposition Library}
  \acro{FROSch}{Fast and Robust Overlapping Schwarz}
  \acro{FE}{finite element}
  \acro{FEM}{finite element method}
  \acro{ASPIN}{Additive Schwarz Preconditioned Inexact Newton}
  \acro{ASPEN}{Additive Schwarz Preconditioned Exact Newton}
  \acro{RASPEN}{Restricted Additive Schwarz Preconditioned Exact Newton}
  \acro{RASPIN}{Restricted Additive Schwarz Preconditioned Inexact Newton}
  \acro{LDC}{lid\hyp{}driven cavity}
  \acro{NKS}{Newton\hyp{}Krylov\hyp{}Schwarz}
  \acro{GMRES}{Generalized Minimal Residual Method}
  \acro{MsFEM}{Multiscale Finite Elements}
  \acro{RGDSW}{Reduced Dimension Generalized Dryja\hyp{}Smith\hyp{}Widlund}
  \acro{GDSW}{Generalized Dryja\hyp{}Smith\hyp{}Widlund}
  \acro{MKL}{Intel Math Kernel Library}
  \acro{BLAS}{Basic Linear Algebra Subprograms}
  \acro{DD}{domain decomposition}
\end{acronym}

\begin{abstract}
Nonlinear Schwarz methods are a type of nonlinear domain decomposition method used as an alternative to Newton's method for solving discretized nonlinear partial differential equations. In this article, the first parallel implementation of a two-level nonlinear Schwarz {\KH method} leveraging the GDSW-type coarse spaces from the Fast and Robust Overlapping Schwarz (FROSch) framework in Trilinos is presented. This framework supports both additive and hybrid two-level nonlinear Schwarz methods and makes use of modifications to the coarse spaces constructed by FROSch to further enhance the robustness and convergence speed of the methods. Efficiency and excellent parallel performance of the software framework {\KH are} demonstrated by applying it to two challenging nonlinear problems: the two-dimensional lid-driven cavity problem at high Reynolds numbers, and a Neo-Hookean beam deformation problem. The results show that two-level nonlinear Schwarz methods scale exceptionally well up to 9\,000 subdomains and are more robust than standard Newton-Krylov-Schwarz solvers for the considered Navier-Stokes problems with high Reynolds numbers or, respectively, for the nonlinear elasticity problems and large deformations. The new parallel implementation provides a foundation for future research in scalable nonlinear domain decomposition methods and demonstrates the practical viability of nonlinear Schwarz techniques for large-scale simulations.
\end{abstract}

\section{Introduction}
\label{sec:intro}

Nonlinear Schwarz methods are robust solvers for discretized nonlinear partial differential equations, such as nonlinear diffusion, nonlinear elasticity, or the Navier–Stokes equations. Instead of solving a discretized problem
\begin{equation}
	\label{eq:nlsys}
	F(u)=0
\end{equation}
with Newton's method directly, an overlapping domain decomposition method is exploited to construct a nonlinearly left-preconditioned system $\mathcal{F}(u):=G(F(u))=0$. This new nonlinear system is then solved with Newton's method, but it shows a different convergence behavior, due to the nonlinear preconditioner $G$.

Many variants have been introduced in the last two decades, including nonlinear one-level Schwarz methods~\cite{caiNonlinearlyPreconditionedInexact2002, doleanNonlinearPreconditioningHow2016} as well as several two-level approaches~\cite{caiNonlinearAdditiveSchwarz2002, parnlschwarz2, hwangClassParallelTwolevel2007, doleanNonlinearPreconditioningHow2016, heinleinAdditiveHybridNonlinear2020, chaouquiLinearNonlinearSubstructured2022, heinleinAdaptiveNonlinearDomain2023}. These yield slightly different nonlinear operators $\mathcal{F}(u)$, but share common goals: faster nonlinear convergence, greater robustness, reduced time to solution, and improved scalability. The first two aspects have been extensively studied in the literature~\cite{caiNonlinearlyPreconditionedInexact2002, doleanNonlinearPreconditioningHow2016, heinleinAdditiveHybridNonlinear2020}, with one key finding being the necessity of an appropriate coarse space (second level) for complex problems and its impact on overall convergence. The latter two aspects, scalability and time to solution, have only been analyzed in~\cite{parnlschwarz1, parnlschwarz2}, for two-level nonlinear Schwarz methods. These implementations do not take advantage of the flexibility of \ac{GDSW}-type coarse spaces and only investigate scalability up to a few hundred subdomains.

To address this gap, we implemented different variants of nonlinear two-level Schwarz in the \ac{FROSch} framework~\cite{heinleinFROSchFastRobust2020}, a collection of linear Schwarz preconditioners within the Trilinos software environment. We provide parallel results for two challenging problems: one based on the stationary Navier–Stokes equations and the other on hyperelasticity. A key result is that our nonlinear Schwarz implementation outperforms the classical Newton approach, where the linearized systems are solved iteratively using a linear Schwarz preconditioner from the \ac{FROSch} package. In particular, the nonlinear Schwarz method demonstrates better weak scalability, faster time to solution in most cases, and higher robustness with respect to large Reynolds numbers or, in the hyperelasticity case, higher loads. We also compare against several nonlinear one-level Schwarz approaches. Our implementation relies on GDSW-type coarse spaces provided in \ac{FROSch}~\cite{heinleinFROSchFastRobust2020}, but includes new modifications designed to accelerate convergence in the nonlinear setting.

Let us note that efficient parallel implementations of different nonlinear domain decomposition methods have already been provided and discussed in the past. A parallel implementation of the one-level nonlinear Schwarz approach ASPIN~\cite{caiNonlinearlyPreconditionedInexact2002} is provided in PETSc. Nonlinear FETI-DP (Finite Element Tearing and Interconnecting - Dual Primal) and BDDC (Balancing Domain Decomposition by Constraints) have also been implemented in parallel and scaled up to {\KH several hundred thousand} cores~\cite{toward,framework}. Other nonlinear domain decomposition methods have been discussed over the past two decades, but we will refrain from providing a detailed overview here and only mention those variants that have been implemented and tested in parallel, since this is the focus of this article.

The remainder of the article is organized as follows: First, we introduce our test problems in Section \ref{sec:problems}, followed by a detailed introduction of the nonlinear Schwarz methods and the coarse spaces in Section \ref{sec:methods}. In Section \ref{sec:implementation}, we introduce the parallel implementation and its structure. Finally, we present parallel results, modifications in the coarse spaces, and some additional key findings in Section \ref{sec:results}.

\section{Nonlinear problems}
\label{sec:problems}
We consider two different test problems: an incompressible flow and a hyperelasticity problem. More precisely, we focus on a lid-driven cavity problem and a beam deformation problem with a Neo-Hookean material model.
The lid-driven cavity problem models the motion of an incompressible fluid in the unit square $\Omega$ using the stationary, dimensionless two-dimensional Navier-Stokes equations. Find $(u,p)\in([H^1_0(\Omega)]^2,L^2_0(\Omega))$ such that
\begin{equation}
	\label{eq:ldc}
	\begin{alignedat}{2}
		-\frac{1}{Re} \Delta u + (u \cdot \nabla)u + \nabla p & =   0 \quad &  & \text{in }\Omega          \\
		\mathrm{div}(u)                                       & =   0       &  & \text{in }\Omega          \\
		u                                                     & = u_0       &  & \text{on }\partial \Omega \\
		p                                                     & =   0       &  & \text{in }(x,y)=(0,0).
	\end{alignedat}
\end{equation}
where $Re$ denotes the Reynolds number, $p$ the pressure, and $u$ the velocity. We impose $u_0=(0,0)$ on the bottom, left, and right boundaries, and $u_0=(1,0)$ on the lid, i.e.~on all nodes with $y=1$. We enforce $p=0$ at the lower-left corner $(0,0)$ to obtain a unique solution. For the velocity, we set the initial iterate $u^{(0)}=0$ in the interior of $\Omega$ and $u^{(0)}=u_0$ on $\partial\Omega$ for both the nonlinear Schwarz and \ac{NKS} solvers. The pressure is initialized with zero everywhere. Because the initial iterate already satisfies the Dirichlet boundary conditions in \eqref{eq:ldc}, we impose homogeneous Dirichlet conditions in every Newton update on $\partial\Omega$ for the velocity and at $(0,0)$ for the pressure.

To model a simple beam deformation, we choose a rectangular domain $\Omega$ with width \qty{5}{\metre} and height \qty{1}{\metre}. Both ends of the beam are fixed by imposing $g_D=0$ on the short edges, while the long edges are traction-free with $g_N=0$. We model the deformation by solving the elasticity problem: find $u\in [H_0^1(\Omega)]^2$ such that
\begin{equation}
	\label{eq:nonlinelas}
	\begin{alignedat}{2}
		-\mathrm{div}(P(F)) & = f_{vol}\quad &  & \text{in }\Omega,           \\
		u                   & = g_D          &  & \text{on }\partial\Omega_D, \\
		n\cdot P(F)         & = g_N          &  & \text{on }\partial\Omega_N,
	\end{alignedat}
\end{equation}
where $F$ is the deformation gradient, and $\Omega_D$ and $\Omega_N$ represent the Dirichlet and Neumann boundaries, respectively. We use a compressible plane-stress Neo-Hookean material model for which the first Piola-Kirchhoff stress tensor is given by
\begin{equation*}
	P(F) = \frac{E}{(1+\nu)}(F-F^{-T}) + \frac{E\nu}{(1+\nu)(1-2\nu)}\mathrm{ln}(\mathrm{det}(F)F^{-T}).
\end{equation*}
Poisson's ratio is set to $\nu=0.3$ and Young's modulus to $E=210\,\mathrm{GPa}$. A volume force is applied perpendicular to the long edges by setting $f_{vol}=(0,-f_y)$ for some scalar $f_y>0$. As before, we start with $u^{(0)}=0$ for both nonlinear Schwarz and \ac{NKS} and impose homogeneous Dirichlet conditions on $\partial\Omega_D$ in every Newton update.

\section{Nonlinear Schwarz methods}
\label{sec:methods}
In this section, we describe the family of nonlinear Schwarz methods. We begin with the one-level nonlinear Schwarz method introduced in \cite{caiNonlinearlyPreconditionedInexact2002}; it forms the foundation for the two-level approaches developed in \cite{heinleinAdditiveHybridNonlinear2020,heinleinAdaptiveNonlinearDomain2023}, which we describe in Section~\ref{sec:nltwo}.
\subsection{The one-level method}
\label{sec:nlone}
The inexact variant of the one-level nonlinear Schwarz method is commonly known as \ac{ASPIN}, while a restricted exact reformulation called \ac{RASPEN} was first presented in \cite{doleanNonlinearPreconditioningHow2016}.

To describe how \ac{ASPIN} and \ac{ASPEN} form the nonlinear {\KH left-preconditioned} system $\mathcal{F}(u)\coloneqq G(F(u)) = 0$, we begin with a partition of $\Omega$ into a set of non-overlapping subdomains, $\Omega_i$ for $i = 1, \dots, N$ such that $\Omega = \cup_{i=1}^{N}\Omega_i$. By adding $k$ layers of finite elements around each subdomain, we construct overlapping subdomains $\Omega_i'$ for $i = 1, \dots, N$, with an overlap width $\delta = kh$, where $h$ represents the diameter of a finite element. We denote by $V(\Omega_i')$ the local finite element spaces associated with each overlapping subdomain $\Omega_i'$. To move between the local and global finite element spaces, we define projection and restriction operators $P_i:V(\Omega_i')\to V(\Omega)$ and $R_i:V(\Omega)\to V(\Omega_i')$. With their help, we define the local nonlinear corrections
\begin{equation}
	\label{eq:local-corrections}
	T_i(u):V(\Omega)\to V(\Omega_i'),\;i=1,\dots,N,
\end{equation}
as the solutions of the local nonlinear problems
\begin{equation}
	\label{eq:local-problems}
	R_iF(u - P_iT_i(u)) = 0.
\end{equation}
The nonlinear preconditioned problem is constructed by combining the local nonlinear corrections on the global domain as follows:
\begin{equation}
	\label{eq:alt-nlsys-1}
	\mathcal{F}_1(u)\coloneqq\sum_{i=1}^NP_iT_i(u) = 0.
\end{equation}
The subscript is to differentiate this one-level method from two-level methods that will be introduced later.
In \cite{caiNonlinearlyPreconditionedInexact2002}, it is shown that if certain assumptions are made on the continuity of $F(u)$, then each exact solution $u^*$ of the original problem \eqref{eq:nlsys} has a neighborhood in which both problems \eqref{eq:nlsys} and \eqref{eq:alt-nlsys-1} share $u^*$ as a unique solution.

In each iteration of Newton's method, which is usually used to solve Equation~\eqref{eq:alt-nlsys-1}, the tangent of $\mathcal{F}_1(u)$ is required. Application of the chain rule to Equation~\eqref{eq:alt-nlsys-1} results in
\begin{equation}
	\label{eq:alt-tangent-1}
	D\mathcal{F}_1(u) = \sum_{i=1}^{N}P_iDT_i(u) = \sum_{i=1}^{N}P_i(R_iDF(u_i)P_i)^{-1}R_iDF(u_i)
\end{equation}
where $u_i\coloneq u-P_iT_i(u)$; see also \cite{caiNonlinearlyPreconditionedInexact2002}. For the sake of simplification, we introduce
\begin{equation}
	\label{eq:q}
	Q_i(u)\coloneq P_i(R_iDF(u)P_i)^{-1}R_iDF(u)
\end{equation}
in order to write $D\mathcal{F}_1(u) = \sum_{i=1}^NQ_i(u_i)$. The exact evaluation of the tangent leads to the \ac{ASPEN} method, whereas approximating the tangent by substituting $u_i$ with $u$ in Equation \eqref{eq:alt-tangent-1} yields the \ac{ASPIN} method. Notably, the \ac{ASPIN} tangent naturally arises from applying a linear one-level overlapping Schwarz preconditioner to $DF(u)$. Thanks to the structure of the \ac{ASPEN} and \ac{ASPIN} tangents, no linear preconditioning is required when solving
$
	D\mathcal{F}_1(u)\delta u = \mathcal{F}_1(u)
$
with a Krylov subspace solver such as \ac{GMRES} \cite{saadGMRESGeneralizedMinimal1986}.

The computation of the local nonlinear corrections $T_i(u)$ involves solving the $N$ nonlinear problems in \eqref{eq:local-problems} using local applications of Newton's method. To reflect the nested structure of Newton's method, we refer to the local applications as \textit{inner Newton's methods} and the global application as the \textit{outer Newton's method}.

Excessive weighting of local nonlinear corrections within overlapping regions can be avoided by using modified prolongation operators during the recombination of the $T_i(u)$ that satisfy the partition of unity property
\begin{equation}
	\label{eq:pou}
	\sum_{i=1}^{N}\widetilde{P}_iR_i = I.
\end{equation}
This modification, first proposed in \cite{doleanNonlinearPreconditioningHow2016}, results in the \ac{RASPEN} method with the modified nonlinear preconditioned problem
$
	\mathcal{F}_{RAS}(u) \coloneqq {\KH \sum_{i=1}^N}\widetilde{P}_iT_i(u) = 0.
$
We use RASPEN on the first level of all our two-level approaches throughout this article.

\subsection{Two-level methods}
\label{sec:nltwo}
To achieve a condition number bound independent of the number of subdomains in the linear case, a global coarse space is needed. Provided the coarse space spans the nullspace of the differential operator of the underlying linear problem, the condition number of the preconditioned system using a two-level additive Schwarz preconditioner has the upper bound
\begin{equation*}
	\kappa(M_{S2}^{-1}A) \leq C(1 + \frac{H}{\delta});
\end{equation*}
see \cite{toselliDomainDecompositionMethods2005,smithDomainDecompositionParallel2004}. A similar behavior can be observed for nonlinear Schwarz methods, for example in \cite{heinleinNonlinearTwoLevelSchwarz2024}, motivating the need for a second level in the nonlinear context as well. The extension of one-level nonlinear Schwarz methods to their two-level counterparts was proposed, for example, in \cite{caiNonlinearlyPreconditionedInexact2002, doleanNonlinearPreconditioningHow2016, parnlschwarz2}. In \cite{hwangClassParallelTwolevel2007}, the authors use a linear second level, which we do not consider here. In Section \ref{sec:coarse}, we provide a detailed description of the coarse spaces and coarse basis functions we use here, but the concept of how to include a coarse space into nonlinear Schwarz methods is quite general and can be explained considering a general finite-dimensional coarse space $V_0$ spanned by some coarse basis functions, a coarse projection operator $P_0: V_0\to V(\Omega)$, and corresponding coarse restriction operator, which we always construct as $R_0 \coloneqq P_0^T$. Similarly to the local nonlinear corrections \eqref{eq:local-corrections}, we define the nonlinear coarse correction
\begin{equation*}
	T_0(u):V(\Omega)\to V_0
\end{equation*}
as the solution of the nonlinear coarse problem
\begin{equation}
	\label{eq:coarse-problem}
	R_0F(u - P_0T_0(u)) = 0.
\end{equation}
The local and coarse nonlinear corrections are the building blocks for a family of two-level nonlinear Schwarz methods. These arise from constructing the preconditioned problem $\mathcal{F}_X(u)$, where the subscript $X$ is a placeholder for a specific method, using different combinations of local and coarse nonlinear corrections. Analogously to the one-level nonlinear Schwarz method, the two-level method is derived by applying Newton's method to solve $\mathcal{F}_X(u)$. Again, this requires solving the global linear system
\begin{equation}
	\label{eq:global-tangent}
	D\mathcal{F}_X(u)\delta u = \mathcal{F}_X(u)
\end{equation}
in each outer Newton iteration, and the linear operator $D\mathcal{F}_X(u)$ inherits a structure similar to a linear two-level Schwarz preconditioned system. Thus, no linear preconditioning is necessary in the two-level nonlinear Schwarz method.

Here, we restrict ourselves to two variants of two-level nonlinear Schwarz. We refer to \cite{heinleinAdditiveHybridNonlinear2020} for further variants and caution the reader that the notation used there differs slightly from ours. Note that these variants could easily be implemented in our framework.

The first method we consider is a simple additive one for which the nonlinear preconditioned problem is defined by \begin{equation*}
	\mathcal{F}_{A}(u) \coloneqq \sum_{i=0}^{N}P_iT_i(u) = P_0T_0(u) + \mathcal{F}_1(u).
\end{equation*}
We call the resulting additive two-level nonlinear Schwarz method simply the \textit{additive} method. In the second variant, the local and coarse nonlinear corrections are coupled multiplicatively by applying the coarse correction first, followed by the additively coupled local corrections. This results in
\begin{equation*}
	\mathcal{F}_{H}(u) \coloneqq \sum_{i=1}^{N}P_iT_i(u-P_0T_0(u)) + P_0T_0(u).
\end{equation*}
We call this two-level nonlinear Schwarz method the \textit{hybrid} method. For the tangents that arise in both methods, we refer to \cite{heinleinAdditiveHybridNonlinear2020}.

\subsection{Coarse spaces}
\label{sec:coarse}
Building on the general description of two-level nonlinear Schwarz methods in the previous section, we proceed in this section with a description of two variations of a concrete coarse basis that we used in our implementation. In particular, building a coarse space actually means defining and computing the prolongation operator $P_0$ from the coarse to the fine \ac{FE} space. Let us note that each column of $P_0$ in matrix representation can be seen as one basis function of the coarse space discretized using the \ac{FE} on the fine scale.

These coarse space basis functions can be constructed either from a \ac{FE} discretization on a coarse triangulation or from the vertices, edges, and, in three dimensions, faces induced by the \ac{DD}. For complex geometries with irregular meshes or partitions, the latter approach is generally preferred, as it yields a coarse space that is, in a sense, custom-built for the specific setting.

An example of a coarse space that follows this principle is the \ac{GDSW} coarse space. It was developed in \cite{dohrmannDomainDecompositionLess2008} in the context of linear two-level Schwarz preconditioners, with a focus on subdomains that only need to satisfy the regularity conditions of a John domain. In later work, the authors proposed the \ac{RGDSW} coarse space as an alternative to the \ac{GDSW} coarse space. They show that, when applied to common problems, the smaller coarse space reduces the time needed to solve the coarse problem while maintaining the condition number estimate of the preconditioner \cite{dohrmannDesignSmallCoarse2017}.

As described in Section \ref{sec:nltwo}, a sufficient condition for the scalability property of (R)GDSW coarse spaces in a linear setting is that they span the nullspace of the underlying differential operator on subdomains that do not touch the Dirichlet boundary \cite{toselliDomainDecompositionMethods2005}. Clearly, this theory does not extend to nonlinear differential operators, much less to two-level nonlinear Schwarz methods. Nevertheless, its application in this context has shown promising results, e.g., in \cite{heinleinAdaptiveNonlinearDomain2023,heinleinAdditiveHybridNonlinear2020}, and we follow this approach in our new parallel implementation.

We now give a brief overview of two variants of the \ac{RGDSW} coarse space and the \ac{GDSW} coarse space in a two-dimensional setting, based on the description found in \cite{heinleinAdditiveHybridNonlinear2020}. To satisfy the nullspace property, a partition of unity is constructed on $\Omega$ and then multiplied by the nullspace. For simplicity, we focus on the scalar Laplace operator, in which case the partition of unity alone spans the nullspace. The partition of unity is obtained in two steps:

\begin{figure}
	\centering
	\begin{subfigure}{0.4\textwidth}
		\centering
		\scalebox{.7}{
			\begin{tikzpicture}[
		declare function={
				f(\x)=(\x<0.25)*(2*\x) +
				(\x>=0.25)*(\x<0.75)*1/2+
				(\x>=0.75)*(\x<1)*(2*\x-1)+
				(\x>=1)*(\x<1.25)*(-2*\x+3)+
				(\x>=1.25)*(\x<1.75)*1/2+
				(\x>1.75)*(-2*\x+4);
			}
	]
	\pgfplotscolormapaccess[0:1]{0}{Set2-5}
	\def\TEMP{\definecolor{mycolor}{rgb}}
	\expandafter\TEMP\expandafter{\pgfmathresult}
	\begin{axis}[
			width=10cm,
			height=8cm,
			view={35}{25},
			xmin=-0.5, xmax=2.5,
			ymin=-0.5, ymax=2.5,
			zmin=0, zmax=1.2,
			xtick=\empty,
			ytick=\empty,
			ztick=\empty,
			grid=none,
			axis lines=none,
			ticks=none,
			clip=false,
		]


		\foreach \i in {0,0.25,...,2} {
				\addplot3[gray!30] coordinates {(\i,0,0) (\i,2,0)};
				\addplot3[gray!30] coordinates {(0,\i,0) (2,\i,0)};
			}

		\foreach \i in {0,...,2} {
				\addplot3[black,line width=0.5pt] coordinates {(\i,0,0) (\i,2,0)};
				\addplot3[black,line width=0.5pt] coordinates {(0,\i,0) (2,\i,0)};
			}

		\draw[variable=\t,domain=0:2,color=mycolor,line width=0.6pt,] plot (axis cs:\t,1,{f(\t)});
		\draw[variable=\t,domain=0:2,color=mycolor,line width=0.6pt,] plot (axis cs:1,\t,{f(\t)});
		\draw[variable=\t,domain=0:0.5,color=black!80,dashed] plot (axis cs:1,0.25,\t);
		\draw[variable=\t,domain=0:0.5,color=black!80,dashed] plot (axis cs:1,0.75,\t);
		\draw[color=mycolor,line width=0.6pt,] (0,0) -- (2,0);
		\draw[color=mycolor,line width=0.6pt,] (0,0) -- (0,2);
		\draw[color=mycolor,line width=0.6pt,] (2,2) -- (2,0);
		\draw[color=mycolor,line width=0.6pt,] (2,2) -- (0,2);

		\addplot3[only marks, mark=*, mark options={scale=0.8, fill=red}] coordinates {(1,1,0)};
		\node[above] at (axis cs:1,1,0) {$\nu$};

	\end{axis}

\end{tikzpicture}}
		\caption{RGDSW}
		\label{fig:rgdsw-function}
	\end{subfigure}
	\hspace{-1mm}
	\begin{subfigure}{0.4\textwidth}
		\centering
		\scalebox{.7}{
			\begin{tikzpicture}[
		declare function={
				f(\x)=(\x<=1)*(1*\x) +
				(\x>1)*(-\x+2);
			}
	]
	\pgfplotscolormapaccess[0:1]{0}{Set2-5}
	\def\TEMP{\definecolor{mycolor}{rgb}}
	\expandafter\TEMP\expandafter{\pgfmathresult}
	\begin{axis}[
			width=10cm,
			height=8cm,
			view={35}{25},
			xmin=-0.5, xmax=2.5,
			ymin=-0.5, ymax=2.5,
			zmin=0, zmax=1.2,
			xtick=\empty,
			ytick=\empty,
			ztick=\empty,
			grid=none,
			axis lines=none,
			ticks=none,
			clip=false,
		]


		\foreach \i in {0,0.25,...,2} {
				\addplot3[gray!30] coordinates {(\i,0,0) (\i,2,0)};
				\addplot3[gray!30] coordinates {(0,\i,0) (2,\i,0)};
			}

		\foreach \i in {0,...,2} {
				\addplot3[black,line width=0.5pt] coordinates {(\i,0,0) (\i,2,0)};
				\addplot3[black,line width=0.5pt] coordinates {(0,\i,0) (2,\i,0)};
			}

		\draw[variable=\t,domain=0:2,color=mycolor,line width=0.6pt,] plot (axis cs:\t,1,{f(\t)});
		\draw[variable=\t,domain=0:2,color=mycolor,line width=0.6pt,] plot (axis cs:1,\t,{f(\t)});
		\draw[color=mycolor,line width=0.6pt,] (0,0) -- (2,0);
		\draw[color=mycolor,line width=0.6pt,] (0,0) -- (0,2);
		\draw[color=mycolor,line width=0.6pt,] (2,2) -- (2,0);
		\draw[color=mycolor,line width=0.6pt,] (2,2) -- (0,2);

		\addplot3[only marks, mark=*, mark options={scale=0.8, fill=red}] coordinates {(1,1,0)};
		\node[above] at (axis cs:1,1,0) {$\nu$};

	\end{axis}

\end{tikzpicture}}
		\caption{MsFEM}
		\label{fig:msfem-function}
	\end{subfigure}
	\vfill
	\begin{subfigure}{0.4\textwidth}
		\centering
		\scalebox{.7}{
			\begin{tikzpicture}[
		declare function={
				f(\x)=(\x<0.75)*(0) +
				(\x>=0.75)*(\x<1)*(4*\x-3)+
				(\x>=1)*(\x<1.25)*(-4*\x+5)+
				(\x>=1.25)*0;
			}
	]
	\pgfplotscolormapaccess[0:1]{0}{Set2-5}
	\def\TEMP{\definecolor{mycolor}{rgb}}
	\expandafter\TEMP\expandafter{\pgfmathresult}
	\begin{axis}[
			width=10cm,
			height=8cm,
			view={35}{25},
			xmin=-0.5, xmax=2.5,
			ymin=-0.5, ymax=2.5,
			zmin=0, zmax=1.2,
			xtick=\empty,
			ytick=\empty,
			ztick=\empty,
			grid=none,
			axis lines=none,
			ticks=none,
			clip=false,
		]


		\foreach \i in {0,0.25,...,2} {
				\addplot3[gray!30] coordinates {(\i,0,0) (\i,2,0)};
				\addplot3[gray!30] coordinates {(0,\i,0) (2,\i,0)};
			}

		\foreach \i in {0,...,2} {
				\addplot3[black,line width=0.5pt] coordinates {(\i,0,0) (\i,2,0)};
				\addplot3[black,line width=0.5pt] coordinates {(0,\i,0) (2,\i,0)};
			}

		\draw[variable=\t,domain=0:2,color=mycolor,line width=0.6pt,] plot (axis cs:\t,1,{f(\t)});
		\draw[variable=\t,domain=0:2,color=mycolor,line width=0.6pt,] plot (axis cs:1,\t,{f(\t)});
		\draw[color=mycolor,line width=0.6pt,] (0,0) -- (2,0);
		\draw[color=mycolor,line width=0.6pt,] (0,0) -- (0,2);
		\draw[color=mycolor,line width=0.6pt,] (2,2) -- (2,0);
		\draw[color=mycolor,line width=0.6pt,] (2,2) -- (0,2);

		\addplot3[only marks, mark=*, mark options={scale=0.8, fill=red}] coordinates {(1,1,0)};
		\node[above] at (axis cs:1,1,0) {$\nu$};

	\end{axis}

\end{tikzpicture}}
		\caption{GDSW vertex}
		\label{fig:gdsw-vertex-function}
	\end{subfigure}
	\hspace{-1mm}
	\begin{subfigure}{0.4\textwidth}
		\centering
		\scalebox{.7}{
			\begin{tikzpicture}[
		declare function={
				f(\x)=(\x<0.25)*(4*\x) +
				(\x>=0.25)*(\x<=0.75)*1+
				(\x>=0.75)*(\x<=1)*(-4*\x+4);
			}
	]
	\pgfplotscolormapaccess[0:1]{0}{Set2-5}
	\def\TEMP{\definecolor{mycolor}{rgb}}
	\expandafter\TEMP\expandafter{\pgfmathresult}
	\begin{axis}[
			width=10cm,
			height=8cm,
			view={35}{25},
			xmin=-0.5, xmax=2.5,
			ymin=-0.5, ymax=2.5,
			zmin=0, zmax=1.2,
			xtick=\empty,
			ytick=\empty,
			ztick=\empty,
			grid=none,
			axis lines=none,
			ticks=none,
			clip=false,
		]


		\foreach \i in {0,0.25,...,2} {
				\addplot3[gray!30] coordinates {(\i,0,0) (\i,2,0)};
				\addplot3[gray!30] coordinates {(0,\i,0) (2,\i,0)};
			}

		\foreach \i in {0,...,2} {
				\addplot3[black,line width=0.5pt] coordinates {(\i,0,0) (\i,2,0)};
				\addplot3[black,line width=0.5pt] coordinates {(0,\i,0) (2,\i,0)};
			}

		\draw[variable=\t,domain=0:2,color=mycolor,line width=0.6pt,] plot (axis cs:1,\t,{f(\t)});
		\draw[variable=\t,domain=0:1,color=black!80,dashed] plot (axis cs:1,0.25,\t);
		\draw[variable=\t,domain=0:1,color=black!80,dashed] plot (axis cs:1,0.75,\t);
		\draw[color=mycolor,line width=0.6pt,] (0,0) -- (2,0);
		\draw[color=mycolor,line width=0.6pt,] (0,0) -- (0,2);
		\draw[color=mycolor,line width=0.6pt,] (2,2) -- (2,0);
		\draw[color=mycolor,line width=0.6pt,] (2,2) -- (0,2);

		\draw[color=mycolor,line width=0.6pt,] (0,1) -- (2,1);
		\draw[color=red,line width=0.8pt,] (1,0) -- (1,1);
		\draw[color=black,line width=0.4pt,] (1.03,0) -- (1.03,1);
		\draw[color=black,line width=0.4pt,] (0.97,0) -- (0.97,1);
	\end{axis}

\end{tikzpicture}}
		\caption{GDSW edge}
		\label{fig:gdsw-edge-function}
	\end{subfigure}
	\caption{RGDSW, MsFEM and GDSW interface vertex functions for the vertex $\nu$ at the interface between four subdomains. The interface is shown in black, the original mesh in grey and the interface vertex function in green. Figure \ref{fig:gdsw-edge-function} shows a GDSW edge function corresponding to the edge marked in red.}
	\label{fig:basis-functions}
\end{figure}
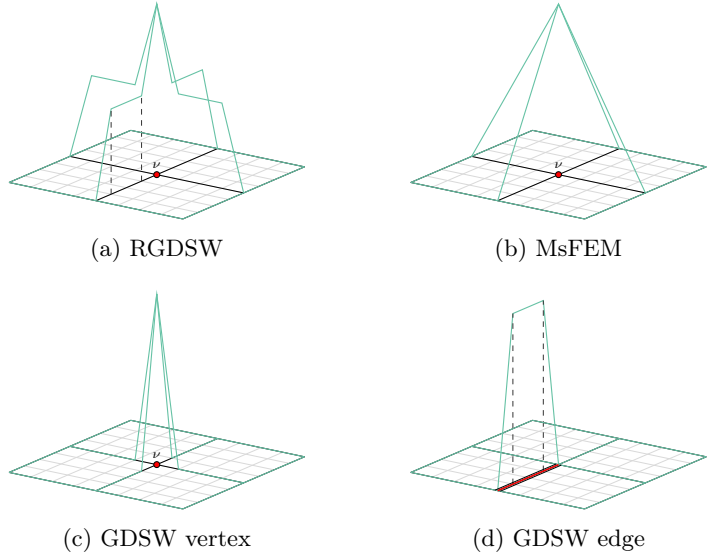

\begin{enumerate}[wide=\parindent, font=\bfseries, label=Step \arabic*:]
	\item We define the \ac{DD} interface
	      $
		      \Gamma \coloneqq \bigcup\limits_{i=1}^{N} \partial\Omega_i
	      $
	      together with the interface disjoint from the true Dirichlet boundary
	      $
		      \Gamma' \coloneqq \Gamma\setminus \partial\Omega_D.
	      $
	      The interface is partitioned into vertices $\mathcal{V}$ and edges $\mathcal{E}$. For the \ac{RGDSW} coarse spaces we define interface vertex functions $\eta_\nu(u):\Gamma\to\mathbb{R}$ on vertices $\nu$ lying in $\Gamma'$ that {\KH fulfill}:
	      \begin{alignat*}{3}
		       &  & \eta_\nu(u)                         & = \delta_{\nu,u}\quad &  & \forall\nu,u\in\mathcal{V},     \\
		       &  & \sum_{\nu\in\mathcal{V}}\eta_\nu(u) & = 1                   &  & \text{on}\ \Gamma',\;\text{and} \\
		       &  & \eta_\nu(u)                         & = 0                   &  & \text{on}\ \partial\Omega_D.
	      \end{alignat*}
	      Here $\delta_{\nu,u}$ denotes the Kronecker delta function. Different variants of the \ac{RGDSW} coarse space result from constructing the interface vertex functions in different ways. In this paper we use two variations that were introduced in \cite{dohrmannDesignSmallCoarse2017} as \textit{Option 1} and \textit{Option 2.2}. We choose to call these variants \ac{RGDSW} and \ac{MsFEM}, respectively, to more closely follow the notation used in \cite{heinleinAdditiveHybridNonlinear2020}. To describe the construction of the interface vertex functions for the \ac{RGDSW} coarse space, we introduce two sets. The first set
	      $
		      \Gamma_\nu \coloneqq \left\{\kappa\in\Gamma\ |\ \exists \epsilon\in\mathcal{E}\ \text{with}\ \kappa,\nu\in\epsilon\right\}
	      $
	      contains all the \ac{FE} nodes in edges connected to the given vertex $\nu\in\mathcal{V}$ including $\nu$ itself. The second set collects all interface vertices within $\Gamma_\nu$, that is,
	      $
		      \mathcal{V}_\nu \coloneqq \mathcal{V} \cap \Gamma_\nu.
	      $
	      Based on these sets the interface vertex functions for the \ac{RGDSW} coarse space are chosen as
	      \begin{equation*}
		      \eta_\nu(u) =
		      \begin{cases}
			      1,           & \text{if}\; u=\nu,                                     \\
			      \frac{1}{2}, & \text{if}\; u\in\Gamma_{\nu}\setminus \mathcal{V}_\nu, \\
			      0,           & \text{otherwise}.
		      \end{cases}
	      \end{equation*}
	      An example of a single interface vertex function for the \ac{RGDSW} coarse space is visualized in Figure \ref{fig:rgdsw-function}. In contrast, the equivalent \ac{MsFEM} interface vertex function is visualized in Figure \ref{fig:msfem-function}. It is constructed by setting
          \begin{equation}\label{eq:inveuclidean}
		      \eta_\nu(u) =
		      \begin{cases}
			      1,                                                 & \text{if}\; u=\nu,                                     \\
			      \frac{1/d_\nu(u)}{1/d_\text{end}(u) + 1/d_\nu(u)}, & \text{if}\; u\in\Gamma_{\nu}\setminus \mathcal{V}_\nu, \\
			      0,                                                 & \text{otherwise},
		      \end{cases}
	      \end{equation}
	      where $d_\nu(u)$ returns the distance between the \ac{FE} node $u$ and the interface vertex $\nu$, and $d_\text{end}(u)$ returns the distance between the node $u$ and the interface vertex at the end of the current interface edge.

	      The \ac{GDSW} coarse space is distinct from the \ac{RGDSW} coarse space in that it defines interface vertex functions as
	      \begin{equation*}
		      \eta_\nu(u) =
		      \begin{cases}
			      1, & \text{if}\; u=\nu, \\
			      0, & \text{otherwise},
		      \end{cases}
	      \end{equation*}
	      and introduces additional interface functions $\eta_\epsilon(u):\Gamma\mapsto\mathbb{R}$ for each edge in $\Gamma'$. To achieve a partition of unity on the interface, the edge functions are chosen to be
	      \begin{equation*}
		      \eta_\epsilon(u) =
		      \begin{cases}
			      1, & \text{if}\; u\in\epsilon\setminus\mathcal{V}, \\
			      0, & \text{otherwise}.
		      \end{cases}
	      \end{equation*}
	      An example of an interface vertex and edge function are shown in Figures \ref{fig:gdsw-vertex-function} and \ref{fig:gdsw-edge-function} respectively.
	\item For simplicity we describe the procedure for constructing \ac{RGDSW}-type coarse basis functions. This extends to \ac{GDSW} coarse basis functions by including the edge functions $\eta_\epsilon$. For each interface vertex $\nu\in\mathcal{V}$, a coarse basis function $\varphi_\nu$ is constructed by extending the corresponding interface vertex function $\eta_\nu$ into the subdomain interiors such that it solves the discrete harmonic problem
	      \begin{equation*}
		      A\varphi_\nu = \begin{pmatrix}A_{II} & A_{I\Gamma} \\ A_{\Gamma I} & A_{\Gamma\Gamma} \end{pmatrix}\begin{pmatrix} \varphi_{I} \\ \varphi_{\Gamma} \end{pmatrix} = 0
	      \end{equation*}
	      where $A$ is the symmetric positive semi-definite matrix that results, for example, from discretizing the Laplace operator with \ac{FEM}. After imposing the boundary conditions, this can be reformulated to
	      \begin{equation*}
		      \begin{pmatrix} \varphi_{I} \\ \varphi_{\Gamma} \end{pmatrix} = \begin{pmatrix}A_{II}^{-1}A_{I\Gamma}\eta_\nu \\ \eta_\nu \end{pmatrix}.
	      \end{equation*}
	      Furthermore, since the interior of the subdomains are independent of each other, $A_{II}$ has a block-diagonal structure:
	      \begin{equation*}
		      A_{II} =
		      \begin{pmatrix}
			      A_{II}^{(1)}                            \\
			       & A_{II}^{(2)}                         \\
			       &              & \ddots                \\
			       &              &        & A_{II}^{(N)}
		      \end{pmatrix}.
	      \end{equation*}
	      This, combined with the local support of $\eta_\nu$, means that
	      \begin{equation*}
		      A_{II}^{-1}A_{I\Gamma}\eta_\nu = \sum_{i\ :\ \nu\in\Omega_i}\bar{R}_i^TA_{II}^{(i)-1}\bar{R}_i A_{I\Gamma}\eta_\nu
	      \end{equation*}
	      with the restriction operator $\bar{R}_i:V(\Omega)\to V(\Omega_i)$. Thus, $A_{II}^{-1}$ can be formed by locally inverting $A_{II}^{(i)}$ on each subdomain adjacent to $\nu$, and the local support property of $\eta_\nu$ carries over to $\varphi_\nu$. This is used in parallel implementations such as \cite{heinleinParallelImplementationTwoLevel2016}. Note that the partition of unity property on interior subdomains results from the $\eta_\nu$ forming a partition of unity on the interface components disjoint from $\partial\Omega$ and the nullspace of the Laplace operator. The coarse projection operator is formed by column-wise concatenation of the $\eta_\nu$.
\end{enumerate}

More generally, the coarse spaces of the type we have described here are constructed by solving the energy minimization problem
\begin{equation}
	\label{eq:min}
	\varphi_\nu \coloneqq \argmin_{v\in V_\nu} |v|^2_{a,\Omega}
\end{equation}
on the set
\begin{equation*}
	V_\nu \coloneqq \{v\in V(\Omega)\ |\ v(u) = \eta_\nu(u)\;\; \forall u\in\Gamma\}
\end{equation*}
with the energy semi-norm
\begin{equation*}
	|v|_{a,\Omega}\coloneqq \sqrt{a(v,v)}
\end{equation*}
induced by the symmetric positive semi-definite bilinear form $a(\cdot,\cdot)$. Since solving the problem \eqref{eq:min} is equivalent to solving
\begin{align*}
	\text{find } \varphi_\nu \in V_\nu \text{ such that } a(\varphi_\nu,v) = 0 \quad \forall v\in V(\Omega)
\end{align*}
the discrete harmonic extension introduced above leads to an energy minimizing extension, and coarse spaces designed in this way satisfy the nullspace property; cf.\ \cite{heinleinAdditiveHybridNonlinear2020}.

\begin{figure}
	\centering
	\begin{subfigure}{0.3\textwidth}
		\centering
		\scalebox{.65}{
			\begin{tikzpicture}[
		declare function={
				f(\x)=(\x<0.25)*(2*\x) +
				(\x>=0.25)*(\x<0.75)*1/2+
				(\x>=0.75)*(\x<1)*(2*\x-1)+
				(\x>=1)*(\x<1.25)*(-2*\x+3)+
				(\x>=1.25)*(\x<1.75)*1/2+
				(\x>1.75)*(-2*\x+4);
			}
	]
	\pgfplotscolormapaccess[0:1]{0}{Set2-5}
	\def\TEMP{\definecolor{mycolor}{rgb}}
	\expandafter\TEMP\expandafter{\pgfmathresult}
	\begin{axis}[
			width=10cm,
			height=8cm,
			view={35}{25},
			xmin=-0.5, xmax=2.5,
			ymin=-0.5, ymax=2.5,
			zmin=0, zmax=1.2,
			xtick=\empty,
			ytick=\empty,
			ztick=\empty,
			grid=none,
			axis lines=none,
			ticks=none,
			clip=false,
		]


		\foreach \i in {0,0.25,...,2} {
				\addplot3[gray!30] coordinates {(\i,0,0) (\i,2,0)};
				\addplot3[gray!30] coordinates {(0,\i,0) (2,\i,0)};
			}

		\foreach \i in {0,...,2} {
				\addplot3[black,line width=0.5pt] coordinates {(\i,0,0) (\i,2,0)};
				\addplot3[black,line width=0.5pt] coordinates {(0,\i,0) (2,\i,0)};
			}

		\draw[variable=\t,domain=0:2,color=mycolor,line width=0.6pt,] plot (axis cs:\t,1,{f(\t)});
		\draw[variable=\t,domain=0:2,color=mycolor,line width=0.6pt,] plot (axis cs:1,\t,{f(\t)});
		\draw[variable=\t,domain=0:0.5,color=black!80,dashed] plot (axis cs:1,0.25,\t);
		\draw[variable=\t,domain=0:0.5,color=black!80,dashed] plot (axis cs:1,0.75,\t);
		\draw[color=mycolor,line width=0.6pt,] (0,0) -- (2,0);
		\draw[color=mycolor,line width=0.6pt,] (0,0) -- (0,2);
		\draw[color=mycolor,line width=0.6pt,] (2,2) -- (2,0);
		\draw[color=mycolor,line width=0.6pt,] (2,2) -- (0,2);

		\addplot3[only marks, mark=*, mark options={scale=0.8, fill=red}] coordinates {(1,1,0)};
		\node[above] at (axis cs:1,1,0) {$\nu$};

	\end{axis}

\end{tikzpicture}}
		\caption{$x$-velocity component}
	\end{subfigure}
	\begin{subfigure}{0.3\textwidth}
		\centering
		\scalebox{.65}{
			\begin{tikzpicture}[
		declare function={
				f(\x)=0;
			}
	]
	\pgfplotscolormapaccess[0:1]{0}{Set2-5}
	\def\TEMP{\definecolor{mycolor}{rgb}}
	\expandafter\TEMP\expandafter{\pgfmathresult}
	\begin{axis}[
			width=10cm,
			height=8cm,
			view={35}{25},
			xmin=-0.5, xmax=2.5,
			ymin=-0.5, ymax=2.5,
			zmin=0, zmax=1.2,
			xtick=\empty,
			ytick=\empty,
			ztick=\empty,
			grid=none,
			axis lines=none,
			ticks=none,
			clip=false,
		]

		\foreach \i in {0,0.25,...,2} {
				\addplot3[gray!30] coordinates {(\i,0,0) (\i,2,0)};
				\addplot3[gray!30] coordinates {(0,\i,0) (2,\i,0)};
			}

		\foreach \i in {0,...,2} {
				\addplot3[black,line width=0.5pt] coordinates {(\i,0,0) (\i,2,0)};
				\addplot3[black,line width=0.5pt] coordinates {(0,\i,0) (2,\i,0)};
			}

		\draw[variable=\t,domain=0:2,color=mycolor,line width=0.6pt,] plot (axis cs:\t,1,{f(\t)});
		\draw[variable=\t,domain=0:2,color=mycolor,line width=0.6pt,] plot (axis cs:1,\t,{f(\t)});
		\draw[color=mycolor,line width=0.6pt,] (0,0) -- (2,0);
		\draw[color=mycolor,line width=0.6pt,] (0,0) -- (0,2);
		\draw[color=mycolor,line width=0.6pt,] (2,2) -- (2,0);
		\draw[color=mycolor,line width=0.6pt,] (2,2) -- (0,2);

		\addplot3[only marks, mark=*, mark options={scale=0.8, fill=mycolor}] coordinates {(1,1,0)};
		\node[above] at (axis cs:1,1,0) {$\nu$};

	\end{axis}

\end{tikzpicture}}
		\caption{$y$-velocity component}
	\end{subfigure}
	\begin{subfigure}{0.3\textwidth}
		\centering
		\scalebox{.65}{
			\begin{tikzpicture}[
		declare function={
				f(\x)=0;
			}
	]
	\pgfplotscolormapaccess[0:1]{0}{Set2-5}
	\def\TEMP{\definecolor{mycolor}{rgb}}
	\expandafter\TEMP\expandafter{\pgfmathresult}
	\begin{axis}[
			width=10cm,
			height=8cm,
			view={35}{25},
			xmin=-0.5, xmax=2.5,
			ymin=-0.5, ymax=2.5,
			zmin=0, zmax=1.2,
			xtick=\empty,
			ytick=\empty,
			ztick=\empty,
			grid=none,
			axis lines=none,
			ticks=none,
			clip=false,
		]

		\foreach \i in {0,0.25,...,2} {
				\addplot3[gray!30] coordinates {(\i,0,0) (\i,2,0)};
				\addplot3[gray!30] coordinates {(0,\i,0) (2,\i,0)};
			}

		\foreach \i in {0,...,2} {
				\addplot3[black,line width=0.5pt] coordinates {(\i,0,0) (\i,2,0)};
				\addplot3[black,line width=0.5pt] coordinates {(0,\i,0) (2,\i,0)};
			}

		\draw[variable=\t,domain=0:2,color=mycolor,line width=0.6pt,] plot (axis cs:\t,1,{f(\t)});
		\draw[variable=\t,domain=0:2,color=mycolor,line width=0.6pt,] plot (axis cs:1,\t,{f(\t)});
		\draw[color=mycolor,line width=0.6pt,] (0,0) -- (2,0);
		\draw[color=mycolor,line width=0.6pt,] (0,0) -- (0,2);
		\draw[color=mycolor,line width=0.6pt,] (2,2) -- (2,0);
		\draw[color=mycolor,line width=0.6pt,] (2,2) -- (0,2);

		\addplot3[only marks, mark=*, mark options={scale=0.8, fill=mycolor}] coordinates {(1,1,0)};
		\node[above] at (axis cs:1,1,0) {$\nu$};

	\end{axis}

\end{tikzpicture}}
		\caption{Pressure component}
	\end{subfigure}
	\caption{Interface values for a monolithic coarse basis function for the velocity component in $x$-direction.}
	\label{fig:monolithic-interface}
\end{figure}
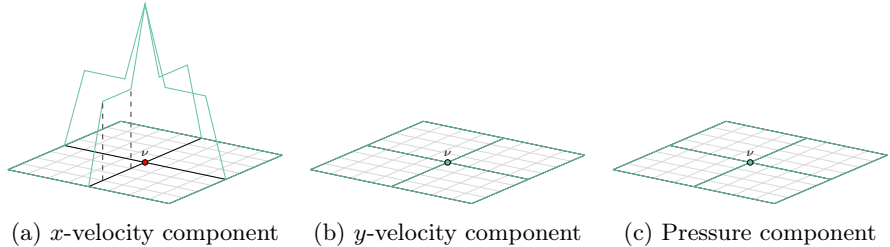

In our two-level nonlinear Schwarz implementation, we use the tangent $DF(u^{(0)})$ of the original problem, evaluated at the initial guess $u^{(0)}$, to compute the extensions into the subdomain interiors. This approach has shown promising results for two-level nonlinear Schwarz methods \cite{heinleinAdditiveHybridNonlinear2020}, but it has not been extensively tested for saddle-point systems so far. For the linearized systems in classical Newton's method, the fully coupled saddle-point system is used to compute the coarse basis functions in \cite{balzaniComputationalFrameworkPharmacomechanical2024, heinleinMonolithicBlockOverlapping2025} for linear two-level Schwarz preconditioners applied to the Navier–Stokes equations. The idea to define the coarse basis functions using the fully coupled saddle point problem was first introduced for Stokes' equations and mixed linear elasticity using standard interpolation instead of the GDSW coarse basis function in \cite{Klawonn:1998:OSM,klawonnComparisonOverlappingSchwarz2000}. The resulting coarse space is commonly referred to as a monolithic coarse space, and we adopt the same approach for our nonlinear two-level Schwarz implementation for the lid-driven cavity problem \eqref{eq:ldc}. Specifically, we construct three complete sets of coarse basis functions: one for the pressure and two for the velocity components, following the procedure described in Section \ref{sec:coarse}. However, each coarse basis function spans all three components. An example of the interface values for an \ac{RGDSW} velocity coarse basis function in the $x$-direction is shown in Figure \ref{fig:monolithic-interface}.

\begin{figure}
	\centering
	\begin{subfigure}{0.3\textwidth}
		\centering
		\includegraphics[width=38mm,height=28mm]{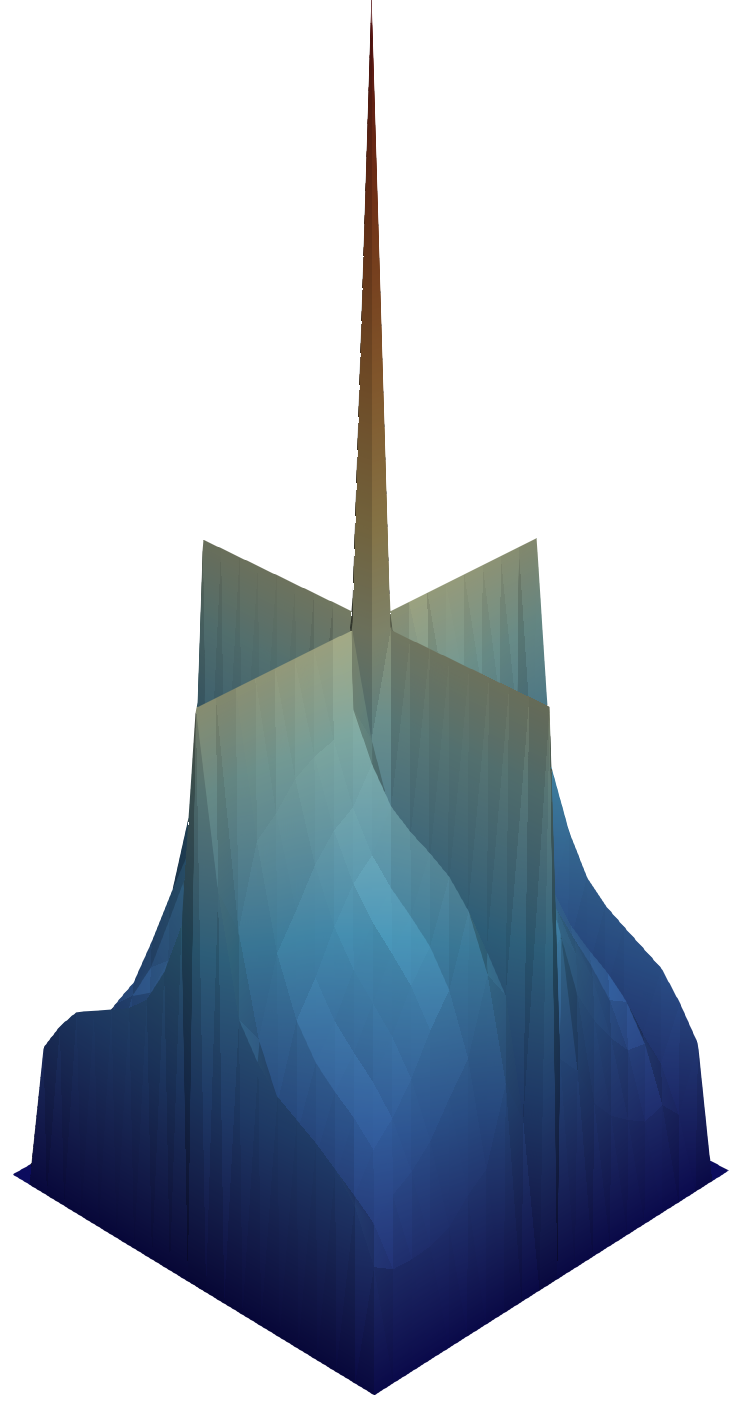}
		\caption{$x$-velocity component}
	\end{subfigure}
	\vspace{-2mm}
	\begin{subfigure}{0.3\textwidth}
		\centering
		\includegraphics[width=38mm,height=28mm]{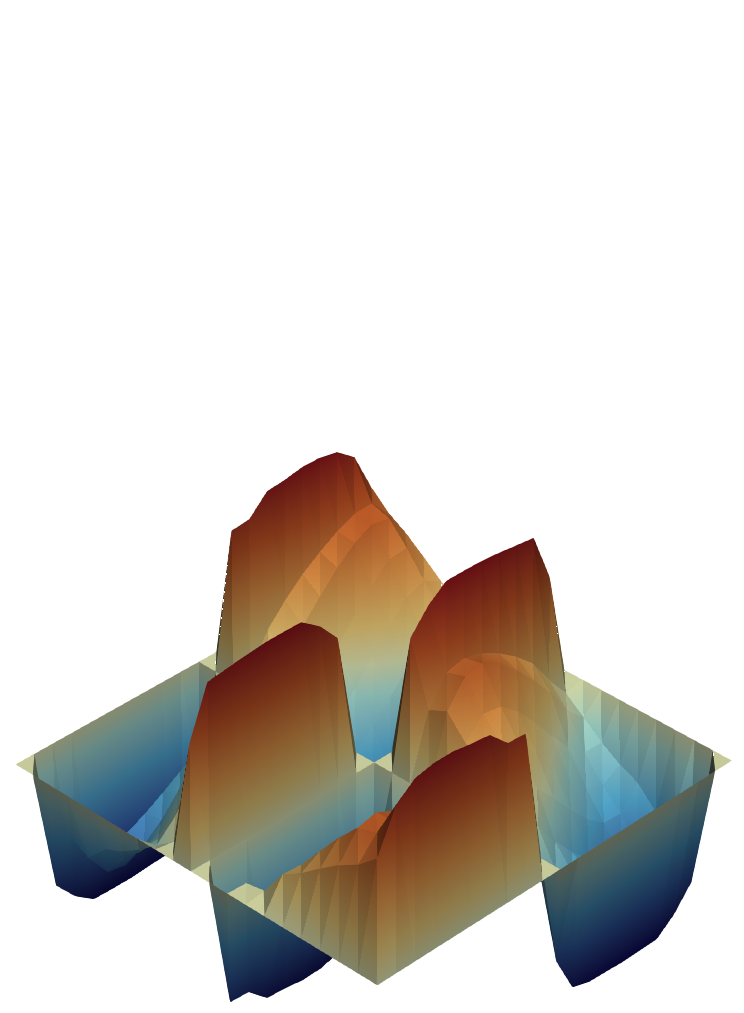}
		\caption{$y$-velocity component}
	\end{subfigure}
	\begin{subfigure}{0.3\textwidth}
		\centering
		\includegraphics[width=38mm,height=28mm]{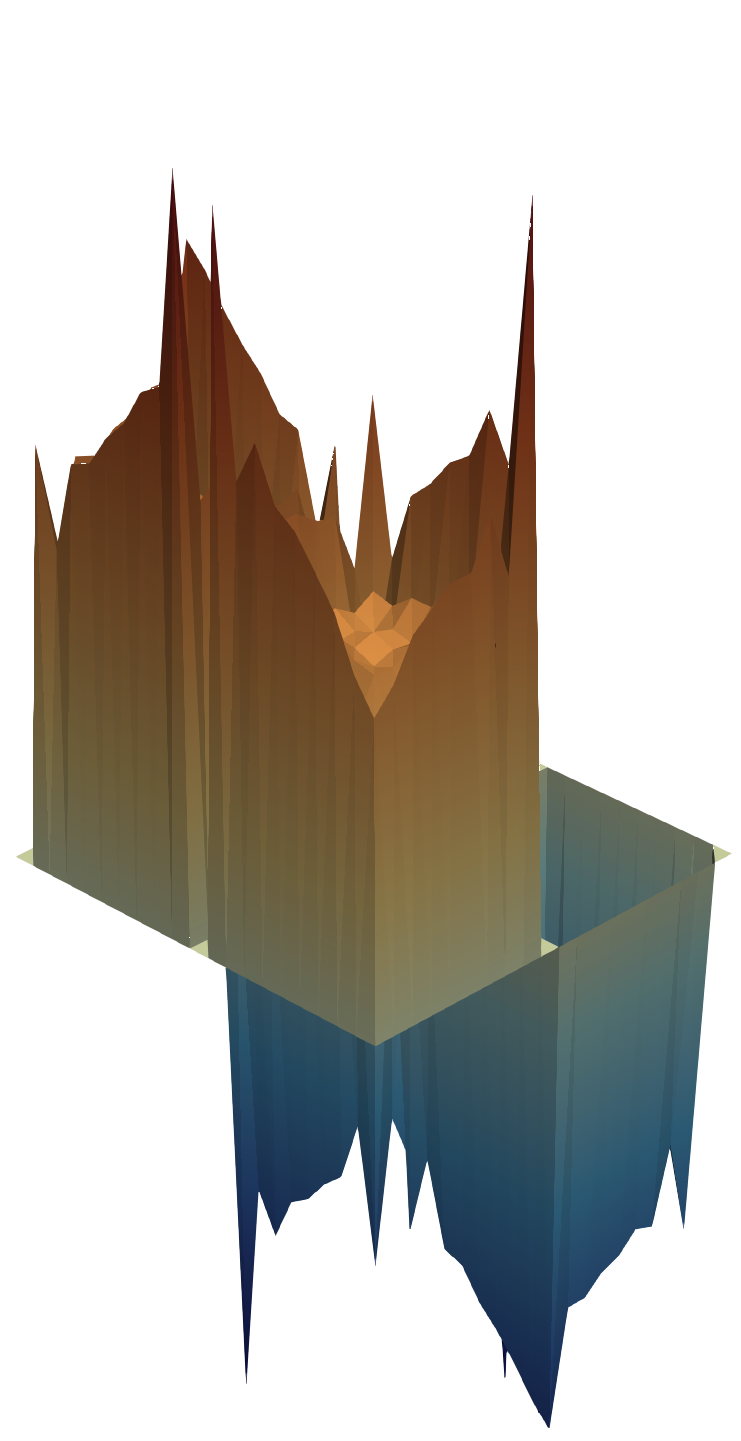}
		\caption{Pressure component}
	\end{subfigure}
	\caption{Complete monolithic coarse basis functions for the velocity component in $x$-direction.}
	\label{fig:monolithic}
\end{figure}

This shows how the partition of unity is constructed for one component at a time, while the other components are set to zero. Extending the interface values into the interior by solving the saddle-point problem results in the coarse basis function shown in Figure \ref{fig:monolithic}. Note that Figure \ref{fig:monolithic} shows the three components of a single coarse basis function, and the full coarse space contains one such basis function for each component at every interface vertex and possibly edge, depending on the coarse space.

Each coarse basis function is represented in the global solution space $V(\Omega)$ by a column vector. Collecting all coarse basis functions column-wise yields a coarse projection operator with the block structure
\begin{equation*}
	P_0 = \begin{pmatrix}
		\varphi_{u,u_0} & \varphi_{u,p_0} \\
		\varphi_{p,u_0} & \varphi_{p,p_0}
	\end{pmatrix}.
\end{equation*}
In our tests, we follow a suggestion made in \cite{balzaniComputationalFrameworkPharmacomechanical2024, heinleinMonolithicBlockOverlapping2025} and set the off-diagonal blocks to zero, as we observed improved performance with this choice. In addition, the coarse space extensions could be recomputed dynamically or at specified iteration intervals using $DF(u^{(k)})$ when the current coarse space no longer performs adequately. Although we did not use this in the present work, we consider it a promising strategy for future problems.

We conclude by noting that other methods for building a coarse space exist. These often require less implementation effort but come with other disadvantages. For example, methods based on linear surrogate models require specialized treatment for each new problem; see \cite{heinleinAdditiveHybridNonlinear2020} where a linear diffusion problem is used to compute the coarse basis functions for nonlinear Schwarz applied to a nonlinear $p$-Laplace problem.

\section{A parallel implementation based on FROSch}
\label{sec:implementation}

In this section, we present our parallel implementation of the two-level nonlinear Schwarz method and its integration into existing software frameworks. Our implementation is based on the \ac{DD} solver package \ac{FROSch} within Trilinos, which provides scalable linear Schwarz preconditioners for large \ac{FE} problems \cite{heinleinFROSchFastRobust2020}. Specifically, we extend the existing operators in \ac{FROSch}, which are implemented using the \mytexttt{Xpetra::Operator} class \cite{trilinos-website}, to support our nonlinear solver.

Most of our implementation and testing is done within a closed-source fork of the \ac{FEDDLib}, an object-oriented \Cpp{} library designed for the development of MPI-parallel \ac{FE} applications~\cite{feddlib}. The \ac{FEDDLib} is tightly integrated with Trilinos, leveraging its distributed linear algebra capabilities, solver infrastructure, and preconditioner packages. In particular, it provides wrappers for Trilinos' parallel linear algebra frameworks via Xpetra, which serves as an abstraction layer for both Epetra and Tpetra. Since Tpetra is the modern Trilinos framework that supports performance portability via Kokkos, our implementation is based primarily on Tpetra to ensure compatibility with current and future high-performance computing architectures. Although we do not yet exploit node-level parallelism, the chosen design readily permits such investigations in the future.

Beyond its linear-algebra interface, the \ac{FEDDLib} offers a framework for \ac{FE} assembly, mesh management, and boundary-condition handling, making it well-suited for testing \ac{DD}-based solvers. Its design ensures efficient handling of overlapping subdomains, which is particularly beneficial for the Schwarz method implemented in \ac{FROSch}. Furthermore, the \ac{FEDDLib} implements a standard nonlinear solver based on Newton's method and supplies interfaces to iterative solvers in Belos \cite{belos} and to preconditioners, most notably \ac{FROSch} via Stratimikos \cite{bartlett2006stratimikos}. Thus, by integrating our solver into the \ac{FEDDLib}, we can apply it to a diverse set of benchmark problems and evaluate its performance in a range of problem settings against existing methods. This integration enables a comprehensive evaluation of robustness, scalability, and efficiency in real-world applications.

In the following, we describe only the new operators we have implemented for the two-level nonlinear Schwarz method and refer to \cite{heinleinParallelImplementationTwoLevel2016,heinleinFROSchFastRobust2020} for descriptions of existing \ac{FROSch} functionality. The paper \cite{heinleinParallelImplementationTwoLevel2016} predates the development of \ac{FROSch} and describes the authors' initial implementation of a two-level linear Schwarz preconditioner based on Trilinos. In particular, it explains the construction of \ac{GDSW} coarse basis functions and how the two-level linear preconditioner is applied through a combination of direct solves on the subdomains and in the coarse space. The construction of the \ac{RGDSW} coarse space basis functions is not explicitly described but follows the same principles as those of the \ac{GDSW} coarse space. The key difference between these coarse spaces is that while \ac{RGDSW} defines a single basis function for each interface vertex, \ac{GDSW} defines separate functions for each vertex, edge, and, in three dimensions, face. This was described in more detail and visualized in Section~\ref{sec:coarse}. In \cite{heinleinFROSchFastRobust2020} \ac{FROSch} is introduced as a restructuring of the earlier implementation, using object-oriented principles to decompose the preconditioner into a set of operator classes derived from \mytexttt{Xpetra::Operator}. This also eases the transition from Epetra to Tpetra. The construction of the coarse space and the use of the preconditioner remain the same after the restructuring. Therefore, we refer the reader to \cite{heinleinFROSchFastRobust2020} for an overview of the existing \ac{FROSch} operators.

Finally, we introduce the Trilinos map, which is the way Trilinos manages distributed linear algebra objects. A map instance stores the global indices of distributed objects, allowing a highly flexible approach to manage distributed matrices and vectors where arbitrary, possibly non-unique, distributions of entries to ranks are possible. We refer to \cite{heinleinFROSchFastRobust2020} for a detailed explanation of the Epetra map class and how a collection of maps is used to manage distributed matrices. This concept also applies directly to Tpetra.

\subsection{Structure of the nonlinear Schwarz implementation}

\begin{figure}[h!]
	\centering
	\begin{tikzpicture}[
    node distance=5mm, 
    every node/.append style={process, font=\footnotesize}, 
    every edge/.style={thick} 
]

\node [fill=orange!10](XpetraOperator) {\texttt{Xpetra::Operator}};
\node [below=of XpetraOperator,fill=orange!10](SchwarzOperator) {\texttt{FROSch::SchwarzOperator}};
\node [right=of SchwarzOperator,thick] (NonLinearOperator) {\texttt{FROSch::NonLinearOperator}};

\draw [arrow] (SchwarzOperator) to (XpetraOperator);

\end{tikzpicture}
	\caption{Abstract classes that form the base of our nonlinear Schwarz implementation.}
	\label{fig:classes-basic}
\end{figure}
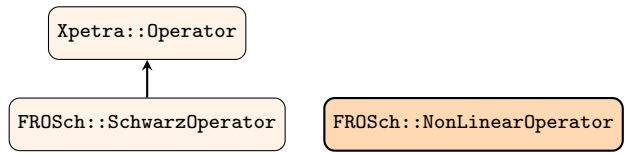

The central building blocks of our implementation are the abstract classes shown in Figure \ref{fig:classes-basic}. The operators implemented for our nonlinear Schwarz solver are highlighted with a thick border and darker shading to distinguish them from those already available in \ac{FROSch}. In the following figures, all operators are within the \ac{FROSch} namespace, so this is not explicitly stated. In \ac{FROSch}, the abstract class \schwarzop{} defines the interface that all concrete operators use. To leverage existing \ac{FROSch} functionality, we follow this architecture. There is one caveat to this design choice: \schwarzop{} is a linear operator interface. The main methods it prescribes are shown in the listing \ref{lst:schwarz}. The \apply{} method is constant, since linear operators can be represented by a constant matrix and applied to any argument. To allow the implementation of nonlinear operators, we have defined the abstract class \nonlinop{}, shown in Listing \ref{lst:nonlinear}. The alternative non-constant \apply{} method is essentially a wrapper for the \texttt{compute()} method, allowing nonlinear operators to be evaluated on arbitrary inputs. By combining these two interfaces, we have a modular structure that implements linear and nonlinear operators with access to the existing \ac{FROSch} functionality.
\newline
\begin{center}
	\begin{minipage}{0.95\linewidth}
		\begin{lstlisting}[language=C++,label=lst:schwarz,caption=Key methods of \schwarzop{}.]
class SchwarzOperator : public Xpetra::Operator {
public:
    virtual int initialize() = 0;
    virtual int compute() = 0;
    virtual void apply(const Xpetra::MultiVector &x,
                       Xpetra::MultiVector &y,
                       ETransp mode=NO_TRANS,
                       double alpha=1.0,
                       double beta=0.0) const;
}
\end{lstlisting}
	\end{minipage}
\end{center}
\begin{center}
	\begin{minipage}{0.95\linewidth}
		\begin{lstlisting}[language=C++,label=lst:nonlinear,caption=The nonlinear operator interface.]
class NonLinearOperator {
public:
    virtual void apply(const Xpetra::MultiVector &x,
                       Xpetra::MultiVector &y,
                       double alpha=1.0,
                       double beta=0.0) = 0;
}
\end{lstlisting}
	\end{minipage}
\end{center}

\subsection{The one-level operators}
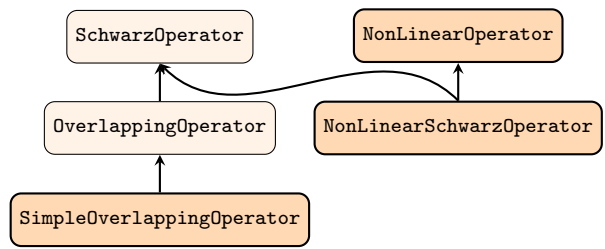
\begin{figure}[h!]
	\centering
	\begin{tikzpicture}[
    node distance=5mm, 
    every node/.append style={process, font=\footnotesize}, 
    every edge/.style={thick} 
]

\node [fill=orange!10](SchwarzOperator) {\schwarzop{}};
\node[below=of SchwarzOperator,fill=orange!10] (OverlappingOperator) {\overlappingop{}};
\node[right=of OverlappingOperator,thick] (NonlinearSchwarzOperator) {\nonlinschwarzop{}};
\node[above=of NonlinearSchwarzOperator,thick] (NonLinearOperator) {\texttt{NonLinearOperator}};
\node[below=of OverlappingOperator,thick] (SimpleOverlappingOperator) {\texttt{SimpleOverlappingOperator}};

\draw [arrow] (OverlappingOperator) to (SchwarzOperator);
\draw [arrow] (NonlinearSchwarzOperator.north) to [out=135,in=315](SchwarzOperator.south);
\draw [arrow] (NonlinearSchwarzOperator) to (NonLinearOperator);
\draw [arrow] (SimpleOverlappingOperator) -- (OverlappingOperator);

\end{tikzpicture}
	\caption{Class hierarchy of our one-level nonlinear Schwarz implementation.}
	\label{fig:classes-one-level}
\end{figure}

We will now discuss the components of our implementation that belong to the one-level nonlinear Schwarz method.

In the outer loop of the nonlinear Schwarz solver, the tangent system \begin{equation}
	\label{eq:tangentsys}
	D\mathcal{F}_1(u^{(k)})\delta u = \mathcal{F}_1(u^{(k)})
\end{equation} must be solved at each iteration. This process requires evaluating the nonlinear operator $\mathcal{F}_1(u)$ and constructing its tangent $D\mathcal{F}_1(u)$ at an arbitrary input $u$. To achieve this, we implemented the \nonlinschwarzop{} for the evaluation and the \simpleoverlappingop{} for the tangent formation. Figure \ref{fig:classes-one-level} shows the inheritance relationship between the new and existing operators.

Recall that $\mathcal{F}_1(u) = \sum_{i=1}^{N}P_iT_i(u)$. Thus, evaluating $\mathcal{F}_1(u)$ requires solving the local nonlinear problems \eqref{eq:local-problems} and combining the local nonlinear corrections $T_i(u)$ on the global domain $\Omega$. To solve the local nonlinear problems, the \nonlinschwarzop{} uses Newton's method with a backtracking line-search globalization as described in \cite{eisenstatChoosingForcingTerms1996}. The tangent systems are solved with a direct solver.

\begin{figure}
	\centering
\resizebox{30mm}{22mm}{%
	\begin{tikzpicture}
		\draw[pattern=crosshatch,pattern color=black] (1,1) -- (4,1) -- (6,3) -- (6,4) -- (7,5) -- (7,6) -- (6,6) -- (6,7) -- (5,7) -- (4,6)  -- (3,6) -- (2,5) -- (2,3) -- (1,2) -- (1,1);
		\draw[fill=green] (2,2) -- (4,2) -- (5,3) -- (5,4) -- (6,5) -- (5,5) -- (5,6) -- (4,5) -- (3,5) -- (3,3);
		\foreach \i [evaluate={\ii=int(\i-1);}] in {0,...,8}{
				\foreach \j [evaluate={\jj=int(\j-1);}] in {0,...,8}{
						\coordinate [shift={(\j,\i)}] (n-\i-\j) at (180:0);
						\ifnum\i>0
							\draw [help lines, semithick, black!70] (n-\i-\j) -- (n-\ii-\j);
						\fi
						\ifnum\j>0
							\draw [help lines, semithick,black!70] (n-\i-\j) -- (n-\i-\jj);
							\ifnum\i>0
								\draw [help lines, semithick,black!70] (n-\i-\j) -- (n-\ii-\jj);
							\fi
						\fi
					}
			}
	\end{tikzpicture}
}
	\vspace{-8mm}
	\caption{The hatched area is the ghost layer corresponding to the green subdomain.}
	\label{fig:ghost-layer}
\end{figure}
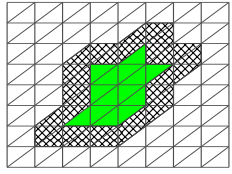

A naive approach to solving $R_iF(u-P_iT_i(u)) = 0$ would involve global assembly of $F(u-P_iT_i(u))$ followed by restriction to each subdomain. However, since the \ac{FE} basis functions have compact support, it is sufficient to perform the assembly locally on each overlapping subdomain $\Omega_i'$, supplemented by an additional layer of \ac{FE} that includes all nodes directly adjacent to $\Omega_i'$. This additional layer, which we denote the \textit{ghost} layer, is visualized in Figure \ref{fig:ghost-layer}. It ensures that all nodes within $\Omega_i'$ are treated as internal nodes during the local assembly of $F(u-P_iT_i(u))$, mirroring their treatment in a global assembly.
To express this formally, we define the set
\begin{equation*}
	\hat{V}_i \coloneqq V(\Omega_i') \oplus \{\textit{FE basis functions on the ghost layer}\}
\end{equation*}
and, assuming the degrees of freedom are ordered correctly, the operator $Z_i: V(\Omega) \to V(\Omega)$ as
\begin{equation*}
	Z_i = \begin{pmatrix} Z_i^{\hat{V}_i} \\ Z_i^{V(\Omega) \setminus \hat{V}_i} \end{pmatrix} \coloneqq \begin{pmatrix} I \\ 0 \end{pmatrix},
\end{equation*}
where $I$ is the identity operator on $\hat{V}_i$ and the zero operator annihilates components outside of $\hat{V}_i$. Due to the compact support of the \ac{FE} basis functions, the following relation holds
\begin{equation*}
	R_i F(u - P_i T_i(u)) = R_i F\big(Z_i(u - P_i T_i(u))\big).
\end{equation*}
Thus, assembling $R_iF(u-P_iT_i(u))$ requires only the part of $u-P_iT_i(u)$ that belongs to $\hat{V}_i$, and once the local assembly is complete, the values corresponding to degrees of freedom in the ghost layer can be discarded. We exploit this fact in our parallel implementation.

In the \ac{FEDDLib}, assembly routines operate in a distributed manner to assemble matrices and vectors according to a unique distribution of mesh elements. This distribution is encoded in a Trilinos map object introduced earlier. To allow \nonlinschwarzop{} to solve nonlinear problems in parallel on subdomains that require the assembly of local tangent matrices, the \apply{} method of \nonlinschwarzop{} calls the assembly routines with serial map objects that are not MPI aware and exist on $\hat{V}_i$. This enables us to use inner Newton methods to obtain solutions
$
	\tilde{u}_i\coloneqq u-P_i(T_i(u))
$
to the local nonlinear problems \eqref{eq:local-problems} on each rank with
$
	\mathrm{supp}(\tilde{u}_i)\in\Omega_i'\oplus\{ghost layer\}.
$
The nonlinear correction is recovered by setting $T_i(u) = R_i(u - \tilde{u})$, i.e.\@, removing the entries of $T_i(u)$ belonging to the ghost layer. Next, the local nonlinear corrections are combined into the distributed preconditioned residual $\mathcal{F}_1(u)$. Our implementation can be used as an \ac{ASPEN} or a \ac{RASPEN}-type method, depending on the user's choice. For \ac{RASPEN}-type methods, we implicitly construct the prolongation operators $\widetilde{P}_i$ that satisfy the partition of unity property \eqref{eq:pou} by averaging the local nonlinear corrections $T_i(u)$ in the overlap. The multiplicity of each degree of freedom required for the averaging, i.e.\@, the number of overlapping subdomains containing that degree of freedom, is determined during initialization of the \nonlinschwarzop{} object.

To evaluate the tangent $D\mathcal{F}_1(u)$ defined in Equation \eqref{eq:alt-tangent-1}, the \simpleoverlappingop{} takes advantage of the work already done by the \nonlinschwarzop{}: the matrices $R_iDF(u-P_iT_i(u))$ are constructed for each overlapping subdomain by the \nonlinschwarzop{} while solving the local nonlinear problems \eqref{eq:local-problems}. When the inner Newton methods have converged, the last tangent is stored for use by the \simpleoverlappingop{}. Again, the compact support of the \ac{FE} basis functions means that it is sufficient to assemble $DF(u-P_iT_i(u))$ only for the degrees of freedom belonging to $\hat{V}_i$ and to restrict the output to $V_i$ by removing entries corresponding to degrees of freedom in the ghost layer. As mentioned in Section \ref{sec:nlone}, the structure of the one-level tangent operator $D\mathcal{F}_1(u)$ is similar to the preconditioned system ${\KH \sum_{i=1}^{N}P_i(R_iDF(u)P_i)^{-1}R_iDF(u)}$ that results from applying a linear Schwarz preconditioner to $DF(u)$. The main difference lies in the replacement of the linear operator $DF(u)$ by shifted counterparts $DF(u - P_iT_i(u))$ on each subdomain. This structural difference prevents the direct use of \texttt{OverlappingOperator}. Instead, the child class \simpleoverlappingop{} overrides the \apply{} function to accommodate the modified structure.

Using the \nonlinschwarzop{} and \simpleoverlappingop{} we are thus able to solve the linear system \eqref{eq:tangentsys} using \ac{GMRES} within the outer Newton method that solves the preconditioned problem \eqref{eq:alt-nlsys-1}.

\subsection{The coarse operators}
\begin{figure}[h!]
	\centering
	\begin{tikzpicture}[
    node distance=5mm, 
    every node/.append style={process, font=\footnotesize}, 
    every edge/.style={thick} 
]

\node [fill=orange!10](SchwarzOperator) {\schwarzop{}};
\node[below=of SchwarzOperator,fill=orange!10] (CoarseOperator) {\coarseop{}};
\node[below=of CoarseOperator,xshift=-2cm,fill=orange!10] (HarmonicCoarseOperator) {\texttt{HarmonicCoarseOperator}};
\node[right=of HarmonicCoarseOperator,thick] (SimpleCoarseOperator) {\simplecoarseop{}};
\node[below=of HarmonicCoarseOperator,fill=orange!10] (IPOUHarmonicCoarseOperator) {\ipouop{}};
\node[right=of IPOUHarmonicCoarseOperator,thick] (NonLinearOperator) {\texttt{NonLinearOperator}};
\node[below=of IPOUHarmonicCoarseOperator,xshift=2cm,thick] (CoarseNonLinearSchwarzOperator) {\texttt{CoarseNonLinearSchwarzOperator}};

\draw [arrow] (CoarseOperator) to (SchwarzOperator);
\draw [arrow] (HarmonicCoarseOperator) to [out=70,in=250](CoarseOperator);
\draw [arrow] (SimpleCoarseOperator.north) to [out=110,in=290](CoarseOperator.south);
\draw [arrow] (IPOUHarmonicCoarseOperator) to (HarmonicCoarseOperator);
\draw [arrow] (CoarseNonLinearSchwarzOperator.north) to [out=110,in=290](IPOUHarmonicCoarseOperator.south);
\draw [arrow] (CoarseNonLinearSchwarzOperator.north) to [out=70,in=250](NonLinearOperator.south);

\end{tikzpicture}
	\caption{Class hierarchy of the coarse operators of our nonlinear Schwarz implementation.}
	\label{fig:classes-two-level}
\end{figure}
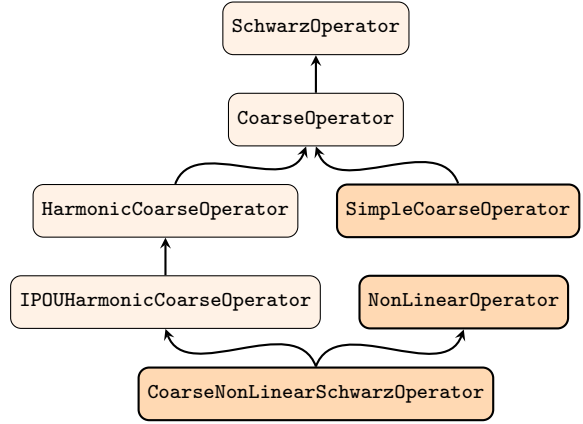

We now proceed to describe the \coarsenonlinschwarzop{} and \simplecoarseop{} that we have implemented to evaluate the nonlinear coarse correction $T_0(u)$. The inheritance structure of these classes is shown in Figure \ref{fig:classes-two-level}. This structure provides \coarsenonlinschwarzop{} with access to \ac{FROSch} functions for constructing the coarse space, applying the resulting prolongation and restriction operators $P_0$ and $R_0$, and inverting the coarse tangent $R_0DF(u-P_0T_0(u))P_0$. The \coarsenonlinschwarzop{} uses Newton's method and its inherited functionality to solve the nonlinear coarse problem \eqref{eq:coarse-problem} for the coarse correction $T_0(u)$. Like the \nonlinschwarzop, it also uses backtracking line-search globalization and a direct solver. Following the approach of the one-level operators, the coarse tangent $R_0DF(u - P_0T_0(u))P_0$ is stored upon convergence of Newton's method and subsequently used to initialize an instance of \simplecoarseop{}. We have implemented \simplecoarseop{} as a wrapper around the existing abstract operator \coarseop{}. This operator implements an \apply{} function that applies the linear operator $Q_0(u_0)$ given $R_0DF(u_0)$ with $Q_0(u)$ from Equation \eqref{eq:q} and $u_0\coloneqq u-P_0T_0(u)$ as before. Thus, by combining \simpleoverlappingop{} and \simplecoarseop{} we can construct the tangent of the preconditioned problem $D\mathcal{F}_X(u)$ for different two-level nonlinear Schwarz variants.

\subsection{Building two-level variants}
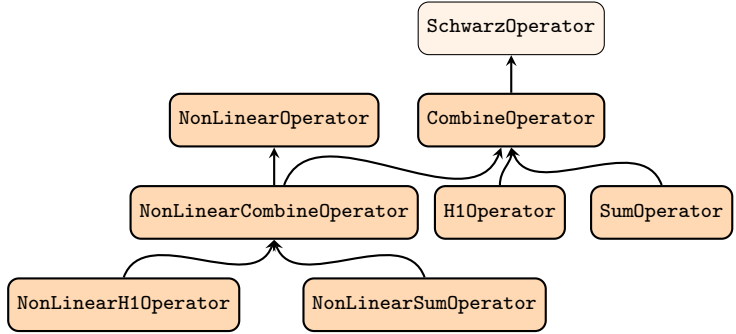
\begin{figure}[h!]
	\centering
	\begin{tikzpicture}[
    node distance=5mm, 
    every node/.append style={process, font=\footnotesize}, 
    every edge/.style={thick} 
]

\node [fill=orange!10](SchwarzOperator) {\schwarzop{}};
\node[below=of SchwarzOperator,thick] (CombineOperator) {\texttt{CombineOperator}};
\node[left=of CombineOperator,thick] (NonLinearOperator) {\texttt{NonLinearOperator}};
\node[below=of NonLinearOperator,thick] (NonLinearCombineOperator) {\texttt{NonLinearCombineOperator}};
\node[right=of NonLinearCombineOperator,xshift=-0.3cm,thick] (H1Operator) {\texttt{H1Operator}};
\node[below=of CombineOperator,xshift=2cm,thick] (SumOperator) {\texttt{SumOperator}};
\node[below=of NonLinearCombineOperator,xshift=-2cm,thick] (NonLinearH1Operator) {\texttt{NonLinearH1Operator}};
\node[below=of NonLinearCombineOperator,xshift=2cm,thick] (NonLinearSumOperator) {\texttt{NonLinearSumOperator}};

\draw [arrow] (CombineOperator) to (SchwarzOperator);
\draw [arrow] (NonLinearCombineOperator) to (NonLinearOperator);
\draw [arrow] (NonLinearCombineOperator) to [out=70,in=250](CombineOperator);
\draw [arrow] (H1Operator.north) to [out=90,in=270](CombineOperator.south);
\draw [arrow] (SumOperator.north) to [out=110,in=290](CombineOperator.south);
\draw [arrow] (NonLinearH1Operator.north) to [out=70,in=250](NonLinearCombineOperator.south);
\draw [arrow] (NonLinearSumOperator.north) to [out=110,in=290](NonLinearCombineOperator.south);

\end{tikzpicture}
	\caption{Class hierarchy of our implementation for combining first and second level operators.}
	\label{fig:classes-variants}
\end{figure}

As described in Section \ref{sec:nltwo}, different two-level nonlinear Schwarz methods arise from the specific construction of the preconditioned problem $\mathcal{F}_X(u) = 0$ from the local and coarse nonlinear corrections $T_i(u),\, i = 0, \dots, N$. To allow a flexible implementation of different nonlinear Schwarz variants, we have designed the operator class structure shown in Figure \ref{fig:classes-variants}. Both \nonlincombineop{} and \combineop{} are abstract operators that define interfaces for evaluating the preconditioned problem $\mathcal{F}_X(u) = 0$ and its tangent, respectively. Specializations of these abstract operators implement two-level additive Schwarz methods. For example, \nonlinsumop{} and \linsumop{} implement the additive method by invoking the one-level and coarse operators and combining their corrections additively. Operators for the hybrid method follow the same approach, and additional operators can be easily implemented to support other combinations of local and coarse corrections.

\subsection{Overlapping mesh partitioning}
\label{subsec:mesh-partitioning}
We will now explain how we create an overlapping partition of a mesh. We assume that we already have a uniquely distributed mesh and only describe how to achieve the overlapping partition. Since \ac{FE} assembly routines are typically implemented to assemble on a per-element basis, it is necessary for each subdomain to know which elements it contains in order to facilitate \ac{FE} assembly on the overlapping subdomains. For this purpose, we build the dual graph of the distributed mesh using ParMETIS \cite{parmetis}. The dual graph is stored as a distributed adjacency matrix. To add a layer of elements to a subdomain, all local column indices are compared to all local row indices. If a column index has no corresponding local row, this indicates that the node in the dual graph represented by that row is not part of the subdomain and must be added. The adjacency matrix is stored as a Trilinos sparse distributed matrix, which simplifies the communication of missing rows. Once the distribution of elements in the overlapping partition of the mesh is known, we communicate the missing nodes in the overlap using Trilinos distributed multivectors and build a map corresponding to the overlapping distribution of nodes. 

In contrast, \ac{FROSch} applies the same technique to the nodal graph of the mesh when building overlapping subdomains. For the linear Schwarz method this is sufficient since the preconditioner is built directly from a pre-assembled distributed matrix. If we followed this approach in our nonlinear Schwarz implementation, it would be necessary to map nodes in the overlap to elements. Without a distributed node to element map, this process does not scale.

\begin{figure}
    \centering
	\begin{subfigure}{0.32\textwidth}
        \centering
        \resizebox{30mm}{22mm}{%
	\begin{tikzpicture}
		\draw[pattern=crosshatch,pattern color=black] (1,1) -- (4,1) -- (6,3) -- (6,4) -- (7,5) -- (7,6) -- (6,6) -- (6,7) -- (5,7) -- (4,6)  -- (3,6) -- (2,5) -- (2,3) -- (1,2) -- (1,1);
		\draw[fill=green] (2,2) -- (4,2) -- (5,3) -- (5,4) -- (6,5) -- (5,5) -- (5,6) -- (4,5) -- (3,5) -- (3,3);
		\foreach \i [evaluate={\ii=int(\i-1);}] in {0,...,8}{
				\foreach \j [evaluate={\jj=int(\j-1);}] in {0,...,8}{
						\coordinate [shift={(\j,\i)}] (n-\i-\j) at (180:0);
						\ifnum\i>0
							\draw [help lines, semithick, black!70] (n-\i-\j) -- (n-\ii-\j);
						\fi
						\ifnum\j>0
							\draw [help lines, semithick,black!70] (n-\i-\j) -- (n-\i-\jj);
							\ifnum\i>0
								\draw [help lines, semithick,black!70] (n-\i-\j) -- (n-\ii-\jj);
							\fi
						\fi
					}
			}
	\end{tikzpicture}
}
        \caption{Nodal: one layer}
	\end{subfigure}
	\begin{subfigure}{0.32\textwidth}
        \centering
        \resizebox{30mm}{22mm}{%
	\begin{tikzpicture}
		\draw[pattern=crosshatch,pattern color=black] (2,1) -- (3,2) -- (3,1) -- (4,2) -- (5,2) -- (5,3) -- (6,4) -- (6,6) -- (4,6) -- (2,4) -- (3,4) -- (2,3) -- (2,1);
		\draw[fill=green] (2,2) -- (4,2) -- (5,3) -- (5,4) -- (6,5) -- (5,5) -- (5,6) -- (4,5) -- (3,5) -- (3,3);
		\foreach \i [evaluate={\ii=int(\i-1);}] in {0,...,8}{
				\foreach \j [evaluate={\jj=int(\j-1);}] in {0,...,8}{
						\coordinate [shift={(\j,\i)}] (n-\i-\j) at (180:0);
						\ifnum\i>0
							\draw [help lines, semithick, black!70] (n-\i-\j) -- (n-\ii-\j);
						\fi
						\ifnum\j>0
							\draw [help lines, semithick,black!70] (n-\i-\j) -- (n-\i-\jj);
							\ifnum\i>0
								\draw [help lines, semithick,black!70] (n-\i-\j) -- (n-\ii-\jj);
							\fi
						\fi
					}
			}
	\end{tikzpicture}
}
        \caption{Dual: one layer}
	\end{subfigure}
	\begin{subfigure}{0.32\textwidth}
        \centering
        \resizebox{30mm}{22mm}{%
	\begin{tikzpicture}
		\draw[pattern=crosshatch,pattern color=black] (1,1) -- (4,1) -- (6,3) -- (6,4) -- (7,5) -- (6,5) -- (7,6) -- (6,6) -- (6,7) -- (5,6) -- (5,7) -- (4,6)  -- (3,6) -- (3,5) -- (2,5) -- (2,3) -- (1,2) -- (2,2) -- (1,1);
		\draw[fill=green] (2,2) -- (4,2) -- (5,3) -- (5,4) -- (6,5) -- (5,5) -- (5,6) -- (4,5) -- (3,5) -- (3,3);
		\foreach \i [evaluate={\ii=int(\i-1);}] in {0,...,8}{
				\foreach \j [evaluate={\jj=int(\j-1);}] in {0,...,8}{
						\coordinate [shift={(\j,\i)}] (n-\i-\j) at (180:0);
						\ifnum\i>0
							\draw [help lines, semithick, black!70] (n-\i-\j) -- (n-\ii-\j);
						\fi
						\ifnum\j>0
							\draw [help lines, semithick,black!70] (n-\i-\j) -- (n-\i-\jj);
							\ifnum\i>0
								\draw [help lines, semithick,black!70] (n-\i-\j) -- (n-\ii-\jj);
							\fi
						\fi
					}
			}
	\end{tikzpicture}
}
        \caption{Dual: two layers}
	\end{subfigure}
    \caption{Different overlaps resulting from adding one layer of \ac{FE} to the subdomain shown in green via the dual graph or the nodal graph.}
	\label{fig:dual-vs-nodal}
\end{figure}
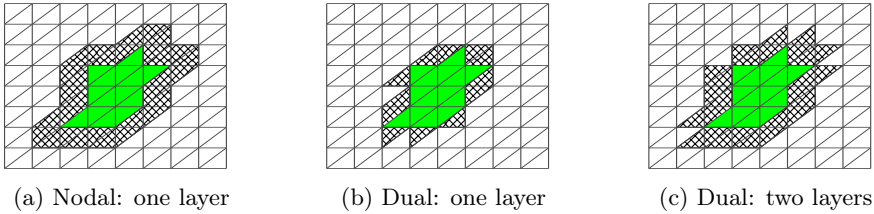

Determining the overlap via the dual or nodal graph results in different overlapping subdomains. This is visualized in Figure \ref{fig:dual-vs-nodal}. Adding twice as many layers via the dual graph compared to the nodal graph, results in similar overlaps. The overlaps differ primarily in the corners of the subdomain. For unstructured meshes this discrepancy might be greater. Nevertheless, we always set the overlap for the nonlinear Schwarz method to be approximately twice as large as for the linear Schwarz method when comparing the relative performance of our nonlinear Schwarz implementation with a standard \ac{NKS} solver.

\subsection{Third-party libraries}
As described before, our implementation uses Trilinos, \ac{FEDDLib}, and ParMETIS. In addition, we use the \ac{MKL} as a \ac{BLAS} backend for Trilinos and its built-in direct solver, PARDISO, for its robust performance. PARDISO is used in sequential mode to solve the linear problem at each Newton iteration within the subdomains and the coarse space \cite{pardiso}. Finally, the HDF5 library is used for parallel writing and storage of the simulation results \cite{The_HDF_Group_Hierarchical_Data_Format}.

\section{Parallel numerical results}
\label{sec:results}
This section begins with an analysis of the weak scalability of our two-level nonlinear Schwarz implementation for the two-dimensional \ac{LDC} problem, followed by a detailed comparison with a standard \ac{NKS} solver as the nonlinearity increases. We then examine the solver's performance on simpler problems and discuss key lessons learned from the tests, particularly the choice of an efficient coarse space and suitable solver tolerances. Our \ac{NKS} implementation uses the \ac{GMRES} algorithm as implemented in the Trilinos package Belos \cite{belos} and a linear Schwarz preconditioner from the aforementioned \ac{FROSch} package. We always use a linear additive two-level Schwarz preconditioner with the same coarse space as the nonlinear Schwarz variant to which we are comparing. We do not use a hybrid linear Schwarz preconditioner because it is not implemented in \ac{FROSch}. However, although the use of an additive linear preconditioner might result in a few more \ac{GMRES} iterations, it does not affect the nonlinear convergence behavior of the \ac{NKS} solver. Furthermore, unless otherwise specified, the \ac{NKS} method uses backtracking line-search globalization as described in \cite{eisenstatChoosingForcingTerms1996}. All applications of backtracking line-search set the tolerance parameter $\bar{\eta}$ and the scaling parameter $t$ to \num{e-3} and the damping parameter $\theta$ to $0.5$. In addition, we enforce an increment tolerance of \num{e-2}, i.e., if the update scaling size $s_k$ drops below this threshold, we stop backtracking. The results presented in this section were obtained using the NHR Fritz supercomputer; see the acknowledgements for more details on Fritz.

\subsection{Lid-driven cavity}
\label{sec:ldc}
The nonlinear Schwarz variants also use the Belos implementation of the \ac{GMRES} algorithm to solve the global tangent system \eqref{eq:global-tangent}. We set the relative tolerance of \ac{GMRES} to \num{e-4}, the maximum number of iterations to $1000$, and the restart value to $500$. For the inner and coarse Newton iterations, we choose a relative tolerance of \num{e-3}, an absolute tolerance of \num{e-14}, and an iteration threshold of a maximum of $10$ iterations. For the outer Newton method, we choose a relative tolerance of \num{e-6} and an absolute tolerance of \num{e-6}. We set the maximum number of iterations to $10$. We discretize the unit square $\Omega$ with a distributed uniform grid. To construct the dual graph of the grid we use ParMETIS \cite{parmetis} as described in the previous section. We use the modified \ac{RGDSW} coarse space described in Section \ref{sec:coarse-observations}, which is constructed to contain a single translation coarse basis function for each degree of freedom at each interface vertex and thus spans the nullspace of the Laplace operator. We found that this coarse space performed much better than any of the other coarse spaces at our disposal; see Section \ref{sec:coarse-observations} for details on tests also comparing the different coarse basis functions. We used an overlap of five for the nonlinear Schwarz solver and chose a subdomain size of $H/h=150$, that is, each subdomain contains at least \num{68000} degrees of freedom excluding the overlap.

\begin{figure}[h!]
	\begin{tikzpicture}
	\pgfplotsset{
		every axis/.append style={
				ybar stacked,
				width=\textwidth,
				height=0.6\textwidth,
				ylabel={Runtime (seconds)},
				xlabel={Num. subdomains},
				symbolic x coords={256, 576, 1024, 2025, 4096, 9216},
				xtick=data,
				enlarge x limits=0.15,
				legend style={at={(1,1)},anchor=north east},
				axis lines*=left, ymajorgrids, yminorgrids,
				ymin=0,
				ymax=250,
				bar width=8pt,
				minor y tick num=1,
				xticklabel style={rotate=0,xshift=0ex,anchor=north},
				cycle list name=Set2-5,
			},
		every axis plot/.append style={
				fill,
			},
	}
  \tikzstyle{mynodestyle} = [rotate=90, anchor=west]

	\begin{axis}[bar shift=-11pt, hide axis]
		\node[mynodestyle] (one) at ([xshift=-11pt]axis cs:256,105) {$450$ $(4)$};
		\node[mynodestyle] at ([xshift=-11pt]axis cs:576,105) {$390$ $(4)$};
		\node[mynodestyle] at ([xshift=-11pt]axis cs:1024,83){$240$ $(3)$};
		\node[mynodestyle] at ([xshift=-11pt]axis cs:2025,85){$210$ $(3)$};
		\node[mynodestyle] at ([xshift=-11pt]axis cs:4096,85){$180$ $(3)$};
		\node[mynodestyle] at ([xshift=-11pt]axis cs:9216,105){$220$ $(4)$};
		\node (oneRe) at ([xshift=-11pt]axis cs:9216,15) {$1$};

        \addplot coordinates  {(256,23.8 ) (576,  22.8 ) (1024, 18.3 ) (2025, 20	) (4096, 17.9)  (9216,22.6  )};   
        \addplot+ coordinates {(256,18.3 ) (576,  19.9 ) (1024, 19.8 ) (2025, 22	) (4096, 25.9)  (9216,31.3)}; 
        \addplot+ coordinates {(256,42.2 ) (576,  40.3 ) (1024, 25 )   (2025, 22.8	) (4096, 20.7)  (9216,25.5)};   
        \addplot+ coordinates {(256,23.7 ) (576,  26   ) (1024, 22.8 ) (2025, 23.9	) (4096, 24.1)  (9216,29.6)};

    \end{axis}

	\begin{axis}[bar shift=0pt, hide axis]

		\node[rotate=0] at (axis cs:256,8) {\scriptsize\color{red}\ding{55}};
		\node[mynodestyle](576) at (axis cs:576,155) {$610$ $(5)$};
		\node[mynodestyle](1024) at (axis cs:1024,120) {$440$ $(4)$};
		\node[mynodestyle](2025) at (axis cs:2025,123) {$370$ $(4)$};
		\node[mynodestyle](4096) at (axis cs:4096,120) {$310$ $(4)$};
		\node[mynodestyle](9216) at (axis cs:9216,115) {$240$ $(4)$};
		\node at (axis cs:9216,15) {$2$};

        \addplot coordinates  {(256, 0)     (576,32.1) (1024, 25.2) (2025, 25.8) (4096, 24.8)  (9216,24  )};
        \addplot+ coordinates {(256, 0)     (576,26.4) (1024, 23.4) (2025, 28.7) (4096, 31.7)  (9216,34.5)};
        \addplot+ coordinates {(256, 0)     (576,66.4) (1024, 46  ) (2025, 40.7) (4096, 36.1)  (9216,29.4  )};
        \addplot+ coordinates {(256, 0)     (576,30.7) (1024, 26.8) (2025, 28.8) (4096, 29.4)  (9216,29.1)};
	\end{axis}

	\begin{axis}[bar shift=11pt]

		\node[xshift=11pt,rotate=0] at (axis cs:256,8) {\scriptsize\color{red}\ding{55}};
		\node[xshift=11pt,rotate=0] at (axis cs:576,8) {\scriptsize\color{red}\ding{55}};
		\node[xshift=11pt,rotate=0] at (axis cs:1024,8) {\scriptsize\color{red}\ding{55}};
		\node [mynodestyle]at([xshift=11pt]axis cs:2025,178) {$620$ $(5)$};
		\node [mynodestyle]at([xshift=11pt]axis cs:4096,135) {$400$ $(4)$};
		\node [mynodestyle]at([xshift=11pt]axis cs:9216,127) {$340$ $(4)$};
		\node(fourRe) at ([xshift=11pt]axis cs:9216,15) {$4$};

        \addplot+ coordinates { (256, 0) (576, 0) (1024, 0) (2025,33.8) (4096,26.6) (9216,24.4)}; 
        \addplot+ coordinates { (256, 0) (576, 0) (1024, 0) (2025,39)  (4096,31.7)  (9216,35.7)};   
        \addplot+ coordinates { (256, 0) (576, 0) (1024, 0) (2025,72.5) (4096,50)   (9216,41.4)};  
        \addplot+ coordinates { (256, 0) (576, 0) (1024, 0) (2025,34.7) (4096,29.7) (9216,28.5)}; 

    \node[anchor=south,text width=1.8cm] (gmres) at (axis cs:256,170){GMRES its. (Newton its.)};

		\legend{
			Inner solve,
			Coarse solve,
			GMRES,
			Other
		}

	\end{axis}

	\node[rotate=0, text width=1.3cm] (Re) at ([yshift=-30,xshift=15]fourRe){Re ($\times 10^{3}$)};

	\draw [thin] (gmres) --  (one);

	\draw [thin] (Re) --  (fourRe);

\end{tikzpicture}
	\caption{Time-to-solution as a function of the number of subdomains for the two-level nonlinear Schwarz solver with Reynolds numbers $Re \in \{1000, 2000, 4000\}$. Runtimes are split into contributions from the inner solver, coarse solver, \ac{GMRES} solver, and remaining overhead. Above each bar, the total number of \ac{GMRES} and Newton iterations are shown as \textit{"GMRES iters. (Newton iters.)".}}
	\label{fig:ldc-schwarz-scalability}
\end{figure}
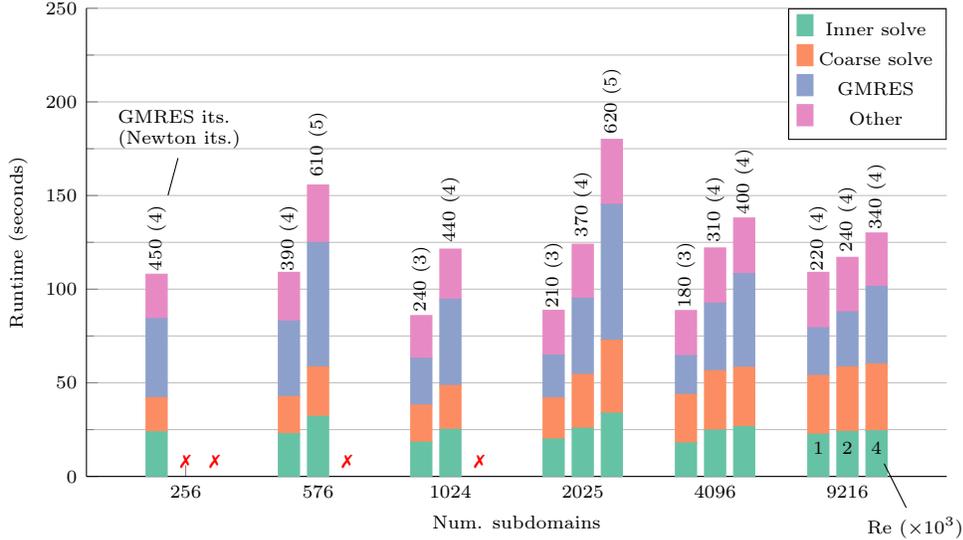

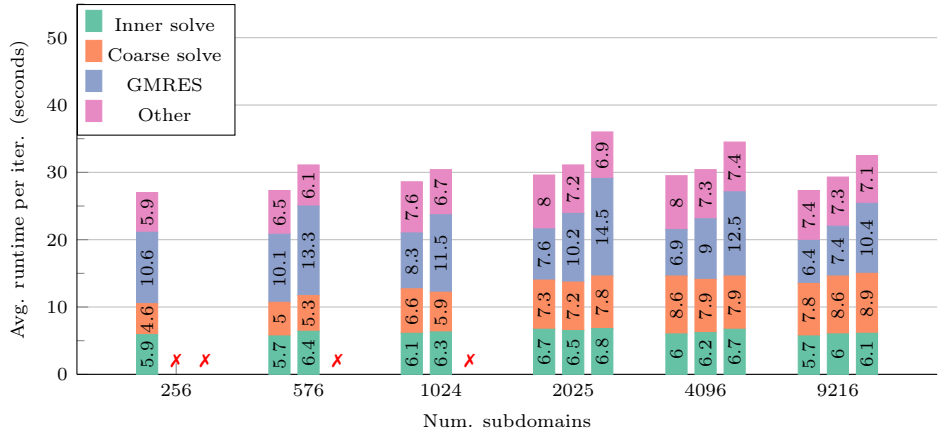
\begin{figure}[h!]
	\begin{tikzpicture}
	\pgfplotsset{
		every axis/.append style={
				ybar stacked,
				width=\textwidth,
				height=0.5\textwidth,
				ylabel={Avg. runtime per iter. (seconds)},
				xlabel={Num. subdomains},
				symbolic x coords={256, 576, 1024, 2025, 4096, 9216},
				xtick=data,
				enlarge x limits=0.15,
				legend style={at={(0,1)},anchor=north west},
				axis lines*=left, ymajorgrids,
				ymin=0,
				ymax=55,
				bar width=8pt,
				nodes near coords,
				nodes near coords style={rotate=90},
				minor y tick num=1,
				xticklabel style={rotate=0,xshift=0ex,anchor=north},
				cycle list name=Set2-5,
			},
		every axis plot/.append style={
				fill,
				every node/.append style={
						text=black,
					},
			},
	}

	\begin{axis}[bar shift=-11pt, hide axis]
		\addplot+ coordinates {(256,5.9  ) (576, 5.7	 ) (1024,   6.1) (2025,  6.7 ) (4096, 6) (9216,5.7 )};
		\addplot+ coordinates {(256,4.6 ) (576, 5 )  (1024,  6.6 ) (2025,  7.3 ) (4096, 8.6 )    (9216,7.8)};
		\addplot+ coordinates {(256,10.6 ) (576, 10.1)  (1024,  8.3 ) (2025,  7.6 ) (4096, 6.9 ) (9216,6.4  )};
		\addplot+ coordinates {(256,5.9 ) (576, 6.5	 ) (1024,   7.6) (2025,  8) (4096, 8   )     (9216,7.4 )};
	\end{axis}

	\begin{axis}[bar shift=0pt, hide axis]

		\node[rotate=0] at (axis cs:256,2) {\scriptsize\color{red}\ding{55}};
		\addplot+ coordinates {(256, 0) (576,6.4 ) (1024, 6.3	) (2025,  6.5   )  (4096,  6.2   ) (9216,6   )};
		\addplot+ coordinates {(256, 0) (576,5.3 ) (1024, 5.9  ) (2025,  7.2  ) (4096,  7.9 )      (9216,8.6)};
		\addplot+ coordinates {(256, 0) (576,13.3) (1024, 11.5  ) (2025,  10.2 ) (4096,  9 )       (9216,7.4  )};
		\addplot+ coordinates {(256, 0) (576,6.1 ) (1024, 6.7	) (2025,  7.2	 )  (4096,  7.3  ) (9216,7.3)};
	\end{axis}

	\begin{axis}[bar shift=11pt]

		\node[xshift=11pt,rotate=0] at (axis cs:256,2) {\scriptsize\color{red}\ding{55}};
		\node[xshift=11pt,rotate=0] at (axis cs:576,2) {\scriptsize\color{red}\ding{55}};
		\node[xshift=11pt,rotate=0] at (axis cs:1024,2) {\scriptsize\color{red}\ding{55}};

		\addplot+ coordinates {(256,0) (576,0) (1024,0)  (2025,  6.8) (4096, 6.7  )   (9216,6.1 )};
		\addplot+ coordinates {(256,0) (576,0) (1024,0)  (2025,  7.8 ) (4096, 7.9)    (9216,8.9  )};
		\addplot+ coordinates {(256,0) (576,0) (1024,0)  (2025,  14.5) (4096, 12.5  ) (9216,10.4)};
		\addplot+ coordinates {(256,0) (576,0) (1024,0)  (2025,  6.9) (4096, 7.4)     (9216,7.1)};

		\legend{
			Inner solve,
			Coarse solve,
			GMRES,
			Other
		}
	\end{axis}
\end{tikzpicture}
	\caption{Time-to-solution averaged over outer Newton iterations. The runtime of each component of the total runtime is shown in each section of the stacked bar chart.}
	\label{fig:ldc-schwarz-scalability-iter}
\end{figure}

Figure \ref{fig:ldc-schwarz-scalability} shows the weak scaling behavior of the two-level nonlinear Schwarz solver when applied to the \ac{LDC} problem. We propose that the reduction in runtime of the method observed for an increasing number of subdomains is due to the growing coarse space, which provides an increasingly good discretization of the original problem. This, in a sense, \textit{shifts} work from the outer Newton method and \ac{GMRES} to the coarse problem, resulting in a better nonlinear preconditioning effect and a decrease in the required outer Newton iterations. Of course, the growing coarse space results in a more difficult and larger coarse problem that takes longer to solve. Figure \ref{fig:ldc-schwarz-scalability-iter} shows the average runtime per outer Newton iteration of the various components of the solver. This makes it easier to see the effects of scaling on the local, coarse, and \ac{GMRES} runtimes. The figure shows that the runtimes of the inner solves remain relatively constant, confirming their perfect weak scalability. The runtime of the coarse solve increases slightly as the coarse problem size increases. This would eventually pose problems for the scalability of the method, in which case alternative strategies would need to be considered for solving the coarse problem, for example, using inexact solvers for the coarse space or three-level nonlinear Schwarz approaches. The \ac{GMRES} runtimes decrease as the linear preconditioning effect of the nonlinear preconditioner improves, which is further verified by the decrease in \ac{GMRES} iterations.

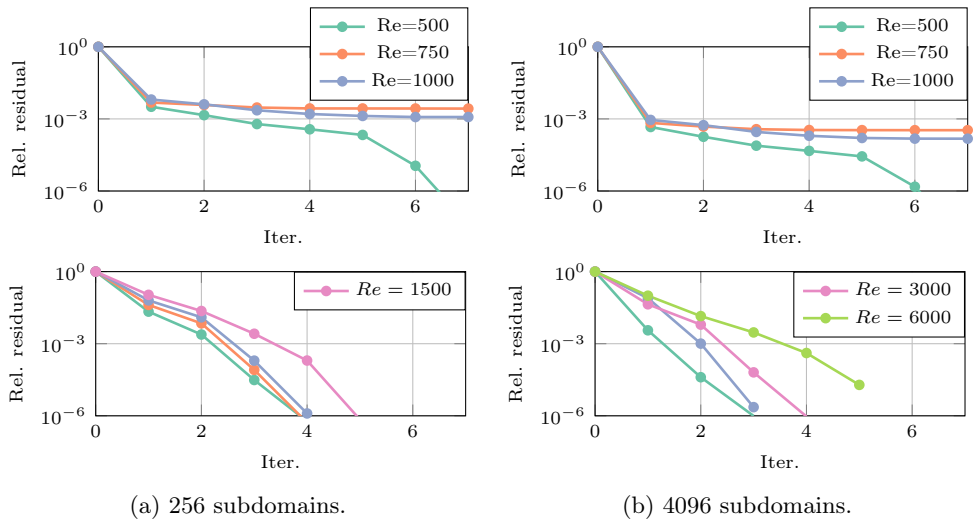
\begin{figure}
	\begin{subfigure}{0.49\textwidth}
		\centering
		\begin{tikzpicture}
    \begin{axis}[
        width=6.5cm, height=3.5cm,
        xlabel={Iter.},
        ylabel={Rel. residual},
        ymode=log,
        grid=major,
        log basis y=10,
        legend style={at={(1.0,0.95)},anchor=east},
				xmin=0,
				xmax=7,
				ymin=1e-6,
				ymax=1,
        every axis plot/.append style={line width=1pt, mark size=1.5pt, mark=*},
    ]
    \addplot coordinates {
        (0, 1) (1, 0.00320133) (2, 0.00143493) (3, 0.000610389) 
        (4, 0.000369944) (5, 0.000215426) (6, 1.11767e-05) (7, 4.18693e-08)
    };
    \addlegendentry{Re=500}
    \addplot coordinates {
        (0, 1) (1, 0.00480191) (2, 0.00380737) (3, 0.00297334)
        (4, 0.00272655) (5, 0.00272013) (6, 0.00270843) (7, 0.00268613)
        (8, 0.00267774) (9, 0.00266959) (10, 0.00243689) (11, 0.0024165)
        (12, 0.00223072) (13, 0.00221035) (14, 0.002201) (15, 0.00219758)
        (16, 0.00219347) (17, 0.00219139) (18, 0.00219044) (19, 0.00219027)
        (20, 0.00218816) (21, 0.00218768) (22, 0.0021868) (23, 0.00218531)
        (24, 0.00218251) (25, 0.00217408) (26, 0.00217307) (27, 0.00217296)
        (28, 0.00217293) (29, 0.00217293)
    };
    \addlegendentry{Re=750}
    \addplot coordinates {
        (0, 1) (1, 0.00640242) (2, 0.00404174) (3, 0.00228814)
        (4, 0.00160117) (5, 0.00132757) (6, 0.00119772) (7, 0.00119687)
        (8, 0.00119643) (9, 0.00119643) (10, 0.00119643) (11, 0.00119643)
        (12, 0.00119643) (13, 0.00119643) (14, 0.00119643) (15, 0.00119643)
        (16, 0.00119643) (17, 0.00119643) (18, 0.00119643) (19, 0.00119643)
        (20, 0.00119643) (21, 0.00119643) (22, 0.00119643) (23, 0.00119643)
        (24, 0.00119643) (25, 0.00119643) (26, 0.00119643) (27, 0.00119643)
        (28, 0.00119643) (29, 0.00119643)
    };
    \addlegendentry{Re=1000}
    \end{axis}
\end{tikzpicture}
		\begin{tikzpicture}
    \begin{axis}[
        group style={group size=1 by 2, vertical sep=1.5cm},
        width=6.5cm, height=3.5cm,
        xlabel={Iter.},
        ylabel={Rel. residual},
        ymode=log,
        grid=major,
        log basis y=10,
        legend style={at={(1.0,1.0)},anchor=north east},
				xmin=0,
				xmax=7,
				ymin=1e-6,
				ymax=1,
        every axis plot/.append style={line width=1pt, mark size=1.5pt, mark=*},
    ]

   \addplot coordinates {
        (0, 1) (1, 0.0216039) (2, 0.00240872) (3, 3.08805e-05) (4, 5.87615e-07)
    };

    \addplot coordinates {
        (0, 1) (1, 0.0410675) (2, 0.00701551) (3, 8.12889e-05) (4, 5.1786e-07)
    };

    \addplot coordinates {
        (0, 1) (1, 0.0635524) (2, 0.0125035) (3, 0.000197537) (4, 1.25259e-06)
    };

    \addplot coordinates {
        (0, 1) (1, 0.10807) (2, 0.0229372) (3, 0.0025968) (4, 0.00019946) (5, 7.0061e-07)
    };
    
    \legend{,,,$Re=1500$}
    \end{axis}
\end{tikzpicture}
		\caption{\num{256} subdomains.}
		\label{fig:residual-ldc-256}
	\end{subfigure}
	\hfill
	\begin{subfigure}{.49\textwidth}
		\centering
		\begin{tikzpicture}
	\begin{axis}[
			group style={group size=1 by 2, vertical sep=1.5cm},
			width=6.5cm, height=3.5cm,
			xlabel={Iter.},
			ylabel={Rel. residual},
			ymode=log,
			grid=major,
			log basis y=10,
			legend style={at={(1.0,0.95)},anchor=east},
			xmin=0,
			xmax=7,
			ymin=1e-6,
			ymax=1,
			every axis plot/.append style={line width=1pt, mark size=1.5pt, mark=*},
		]
		\addplot coordinates {
				(0, 1) (1, 0.000461063) (2, 0.000180623) (3, 7.69231e-05) (4, 4.67144e-05) (5, 2.75539e-05) (6, 1.50771e-06) (7, 6.3394e-09)
			};
		\addlegendentry{Re=500}
		\addplot coordinates {
				(0, 1) (1, 0.000683264) (2, 0.000480095) (3, 0.00037343) (4, 0.000345118) (5, 0.000341495) (6, 0.000339957) (7, 0.000339534) (8, 0.000414796) (9, 0.000364827) (10, 0.000333536) (11, 0.000317937) (12, 0.000275248) (13, 0.000256005) (14, 0.000254507) (15, 0.000254237) (16, 0.000262081) (17, 0.00025798) (18, 0.000261379) (19, 0.000260608)
			};
		\addlegendentry{Re=750}
		\addplot coordinates {
				(0, 1) (1, 0.000907096) (2, 0.000549946) (3, 0.000287238) (4, 0.000198314) (5, 0.000160798) (6, 0.000151072) (7, 0.000150393) (8, 0.000150489) (9, 0.000155954) (10, 0.000154682) (11, 0.000154337) (12, 0.000159684) (13, 0.000159413) (14, 0.000172403) (15, 0.000172074) (16, 0.000172708) (17, 0.000189149) (18, 0.000190431) (19, 0.000198664)
			};
		\addlegendentry{Re=1000}
	\end{axis}
\end{tikzpicture}
		\begin{tikzpicture}
	\begin{axis}[
			group style={group size=1 by 2, vertical sep=1.5cm},
			width=6.5cm, height=3.5cm,
			xlabel={Iter.},
			ylabel={Rel. residual},
			ymode=log,
			grid=major,
			log basis y=10,
			legend style={at={(1.0,1.0)},anchor=north east},
			xmin=0,
			xmax=7,
			ymin=1e-6,
			ymax=1,
			every axis plot/.append style={line width=1pt, mark size=1.5pt, mark=*},
		]
		\addplot coordinates {
				(0, 1) (1, 0.00361003) (2, 4.0282e-05) (3, 9.87625e-07)
			};

		\pgfplotsset{cycle list shift=1}
		\addplot coordinates {
				(0, 1) (1, 0.0783692) (2, 0.00100882) (3, 2.29151e-06)
			};

		\pgfplotsset{cycle list shift=1}
		\addplot coordinates {
				(0, 1) (1, 0.0436383) (2, 0.00622051) (3, 6.42127e-05) (4, 9.38152e-07)
			};

		\addplot coordinates {
				(0, 1) (1, 0.0994514) (2, 0.0140707) (3, 0.00294689) (4, 0.000409486) (5, 1.93227e-05)
			};
        \legend{,,$Re=3000$,$Re=6000$}
	\end{axis}
\end{tikzpicture}
		\caption{\num{4096} subdomains.}
		\label{fig:residual-ldc-4096}
	\end{subfigure}
	\caption{Convergence behavior of \ac{NKS} (top row) and nonlinear two-level Schwarz (bottom row) for various Reynolds numbers.}
	\label{fig:residual-ldc}
\end{figure}

To see how well our two-level nonlinear Schwarz implementation handles nonlinearities in comparison with a standard \ac{NKS} approach, we set the solver parameters as described for the weak scalability test and choose a problem size of $256$ subdomains with $H/h=150$ as before. We set the overlap of the \ac{NKS} solver to two instead of five. This discrepancy results in a similar physical overlap as described in Section \ref{subsec:mesh-partitioning}. We solve the \ac{LDC} problem and incrementally increase the Reynolds number. Figure \ref{fig:residual-ldc-256} shows the relative residuals of the \ac{NKS} method and the nonlinear Schwarz method respectively for a number of different Reynolds numbers. For $Re= 500$ the \ac{NKS} method converges in seven iterations, but for $Re = 750$ it plateaus after a few iterations and fails to meet the absolute and relative tolerances of \num{e-6}. In contrast, the nonlinear Schwarz method converges in four iterations for $Re= 500$ and continues to converge up to $Re = 1500$. As noted previously in this section, our two-level nonlinear Schwarz implementation generally performs better with more subdomains. This observation is verified in Figure \ref{fig:residual-ldc-4096} that shows the relative residuals that result from repeating the previous test with \num{4096} subdomains. As expected, \ac{NKS} displays unchanged nonlinear {\KH behavior}, whereas the nonlinear Schwarz method solves the lid-driven cavity problem up to $Re=6000$.

To conclude this section, we would like to demonstrate the utility of two-level nonlinear Schwarz methods compared to the one-level method. It was noted in \cite{parnlschwarz2} that a second level is required to achieve a nonlinear Schwarz solver that scales in both a linear and nonlinear sense. We verified this for our implementation on a nonlinear diffusion problem in \cite{heinleinNonlinearTwoLevelSchwarz2024}. Furthermore, in this section we demonstrated that when applied to the \ac{LDC} problem our two-level nonlinear Schwarz implementation not only scales, but even improves in both its linear and nonlinear performance with an increasing number of subdomains. In contrast to this, Table \ref{tab:aspen} shows the inner nonlinear and outer linear and nonlinear iterations required by the one-level method. For \num{16} subdomains the one-level method fails to converge for $Re>100$. When the run is scaled in a weak sense to \num{64} subdomains, linear preconditioning deteriorates; Table~\ref{tab:aspen-64} shows the resulting jump in linear iterations. For \num{256} subdomains the \ac{GMRES} solver fails to converge after \num{1000} iterations even for $Re=10$. Thus, we have verified also for the \ac{LDC} problem that we require a two-level method to achieve scalability in the outer linear solver. We argue that the nonlinear scalability of the one-level method, as claimed in \cite{parnlschwarz2}, should be tested for more subdomains, however, this is not possible since the linear solver breaks down so rapidly.

\begin{table}[htbp]
	\centering
	\begin{subtable}[t]{0.45\textwidth}
		\centering
		\renewcommand{\arraystretch}{1.2} 
		{\small\begin{tabular}{|r|r|r|r|}
				\hline
				\multicolumn{1}{|l|}{$Re$} &
				\multicolumn{1}{l|}{Outer} &
				\multicolumn{1}{l|}{Inner} &
				\multicolumn{1}{l|}{GMRES}                 \\
				\hline
				10                         & 3 & 2.2 & 260 \\ \hline
				100                        & 4 & 4.2 & 550 \\ \hline
				500                        & - & -   & -   \\ \hline
			\end{tabular}}
		\caption{\num{16} subdomains.}
		\label{tab:aspen-16}
	\end{subtable}
	\hspace{1mm}
	\begin{subtable}[t]{0.35\textwidth}
		\centering
		\renewcommand{\arraystretch}{1.2} 
		{\small\begin{tabular}{|r|r|r|}
				\hline
				\multicolumn{1}{|l|}{Outer} &
				\multicolumn{1}{l|}{Inner}  &
				\multicolumn{1}{l|}{GMRES}               \\
				\hline
				3                           & 2.1 & 800  \\ \hline
				4                           & 3.5 & 4100 \\ \hline
				-                           & -   & -    \\ \hline
			\end{tabular}}
		\caption{\num{64} subdomains.}
		\label{tab:aspen-64}
	\end{subtable}
	\caption{Performance of one-level nonlinear Schwarz (ASPEN) for various Reynolds numbers and subdomain counts. A subdomain size of $H/h = 127$ was used to match the tests in \cite{caiNonlinearlyPreconditionedInexact2002}. The same solver settings were used as for the tests discussed previously in this section. Missing values indicate that the solver diverged.}
	\label{tab:aspen}
\end{table}

Note that for Table \ref{tab:aspen-16} we have chosen the same subdomain size and number of subdomains as used in \cite{caiNonlinearlyPreconditionedInexact2002} to test the \ac{ASPIN} method on the {\KH two-dimensional} \ac{LDC} problem. We obtained different results than in \cite{caiNonlinearlyPreconditionedInexact2002} for high Reynolds numbers{\KH , even when using \ac{ASPEN}, which} performed better than \ac{ASPIN} in our experience. We assume that this observed performance difference stems from the different problem definitions and discretization schemes used: velocity-vorticity with finite differences in \cite{caiNonlinearlyPreconditionedInexact2002}, while we use the standard velocity-pressure formulation with \ac{FE}. Additionally, in~\cite{caiNonlinearlyPreconditionedInexact2002}, a stronger, but also more expensive, globalization strategy was used. For the more frequently used velocity-pressure formulation discretized with finite elements, however, it is clear that a second level is indispensable.

\subsection{Hyperelastic beam}
\label{sec:2dbeam}
We now analyze the performance of the nonlinear Schwarz solver for the nonlinear elasticity problem \eqref{eq:nonlinelas}, designed as a sanity check to assess strong scalability and robustness. The \ac{MsFEM} coarse space is used with the modification from Section \ref{sec:coarse-observations}. In this case, we construct the coarse space to include translation and linearized rotation to span the nullspace of the differential operator. The outer Newton method is set to a relative tolerance of \num{e-4} with a maximum of \num{10} iterations. \ac{GMRES} uses a relative tolerance of \num{e-6}, with the maximum number of iterations set to \num{100} and no restarting. The inner and coarse nonlinear solvers use a relative tolerance of \num{e-3}, an absolute tolerance of \num{e-9}, and a maximum of \num{15} iterations. The overlap is set to \num{10} for the nonlinear Schwarz solver and five for the \ac{NKS} solver, ensuring a similar physical overlap as explained earlier. Line-search backtracking is deactivated in all solvers in this case because we observed improved performance without it.

\begin{figure}[h!]
	\begin{tikzpicture}
	\pgfplotsset{
		every axis/.append style={
				ybar stacked,
				width=\textwidth,
				height=0.5\textwidth,
				ylabel={Runtime (seconds)},
				xlabel={Num. subdomains},
				symbolic x coords={72, 144, 288, 576, 1152},
				xtick=data,
				enlarge x limits=0.15,
				legend style={at={(1,1)},anchor=north east},
				axis lines*=left, ymajorgrids, yminorgrids,
				ymin=0,
				ymax=100,
				bar width=8pt,
				minor y tick num=1,
				xticklabel style={rotate=0,xshift=0ex,anchor=north},
				cycle list name=Set2-5,
			},
		every axis plot/.append style={
				fill,
			},
	}
	\tikzstyle{mynodestyle} = [rotate=90, anchor=west]

	\begin{axis}[bar shift=11pt, hide axis]
		\node[mynodestyle] at ([xshift=11pt]axis cs:72,58) {$65$ $(3)$};
		\node[mynodestyle] (one) at ([xshift=11pt]axis cs:144,27.5) {$62$ $(3)$};
		\node[mynodestyle] at ([xshift=11pt]axis cs:288,13.6){$53$ $(3)$};
		\node[mynodestyle] at ([xshift=11pt]axis cs:576,6.8){$44$ $(3)$};
		\node[mynodestyle] at ([xshift=11pt]axis cs:1152,4.5){$42$ $(3)$};
		\node[mynodestyle](oneRe) at ([xshift=11pt]axis cs:72,1) {H};

		\addplot coordinates {(72,58) (144,27.5) (288,13.6) (576,6.8) (1152,4.5)};
	\end{axis}

	\begin{axis}[bar shift=0pt, hide axis]
		\node[mynodestyle] at (axis cs:72,80) {$100$ $(4)$};
		\node[mynodestyle](two) at (axis cs:144,34.5) {$97$ $(4)$};
		\node[mynodestyle]at (axis cs:288,17) {$81$ $(4)$};
		\node[mynodestyle]at (axis cs:576,8) {$75$ $(4)$};
		\node[mynodestyle]at (axis cs:1152,5) {$71$ $(4)$};
		\node [mynodestyle]at (axis cs:72,1) {A};

		\addplot+ coordinates {(72,81.2) (144,35.6) (288,16.6) (576,8.6) (1152,5.4)};
	\end{axis}


	\begin{axis}[bar shift=-11pt]
		%

		\node [mynodestyle]at([xshift=-11pt]axis cs:72,33) {$162$ $(5)$};
		\node [mynodestyle](three) at([xshift=-11pt]axis cs:144,16) {$159$ $(5)$};
		\node [mynodestyle]at([xshift=-11pt]axis cs:288,8) {$136$ $(5)$};
		\node [mynodestyle]at([xshift=-11pt]axis cs:576,4) {$124$ $(5)$};
		\node [mynodestyle]at([xshift=-11pt]axis cs:1152,3.2) {$108$ $(5)$};
		\node [mynodestyle]at([xshift=-11pt]axis cs:72,1) {NKS};

		\addplot+ coordinates {(72,34) (144,17.1) (288,9.1) (576,5.1) (1152,4.2)};
		\node[text width=1.8cm] (gmres) at (axis cs:144,75){GMRES its. (Newton its.)};

	\end{axis}

	\draw [thin] (gmres) --  (one);
	\draw [thin] (gmres) --  (two);
	\draw [thin] (gmres) --  (three);

\end{tikzpicture}
	\caption{Strong scalability of \ac{NKS} and the two nonlinear Schwarz variants when applied to the nonlinear elasticity problem \eqref{eq:nonlinelas} with a volume force of $f_y = 4$ MN/m$^{2}$.}
	\label{fig:elasticity-strong}
\end{figure}
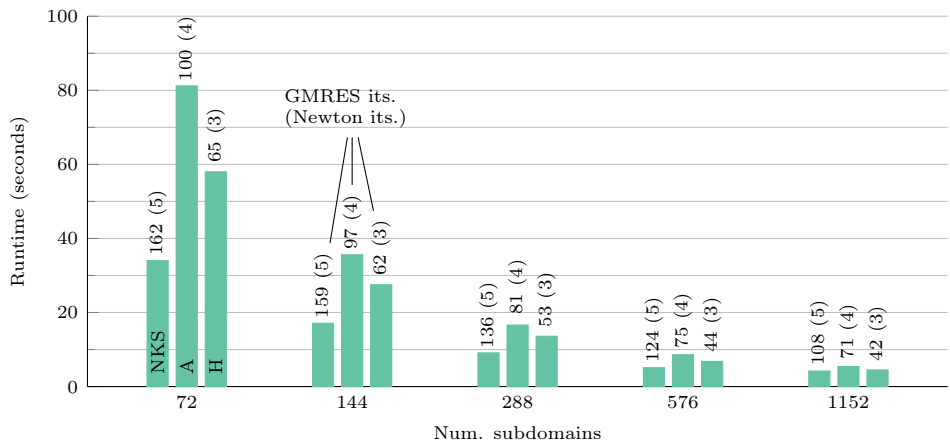

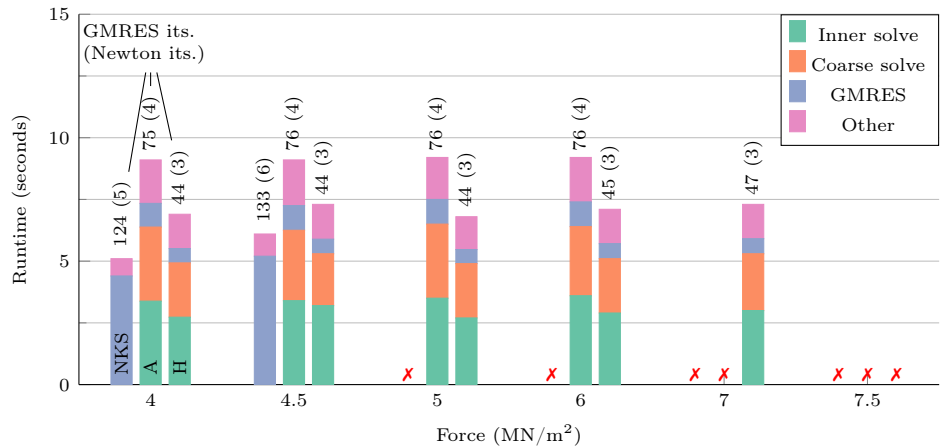
\begin{figure}[h!]
	\begin{tikzpicture}
	\pgfplotsset{
		every axis/.append style={
				legend style={at={(1,1)},anchor=north east},
				axis lines*=left, ymajorgrids, yminorgrids,
				width=\textwidth, height=0.5\textwidth,
				ymin=0,
				ymax=15,
				ybar stacked,
				bar width=8pt,
				minor y tick num=1,
				symbolic x coords={1,2,3,4,5,6},
				xtick={1,2,3,4,5,6},
				xticklabels from table={\hybrid}{Force},
				xticklabel style={rotate=0,xshift=0ex,anchor=north},
				ylabel={Runtime (seconds)},
				xlabel={Force (MN/m$^{2}$)},
				cycle list name=Set2-5,
			},
		every axis plot/.append style={fill},
	}

	\tikzstyle{mynodestyle} = [rotate=90, anchor=west]

	\pgfplotstableread{
		Location Force  GlobalSolve   InnerSolve   CoarseSolve   GMRES    Other
		1        4      6.9           2.73         2.2           0.58     1.39
		2        4.5    7.3           3.2          2.1           0.59     1.41
		3        5      6.8           2.7          2.2           0.57     1.33
		4        6      7.1           2.9          2.2           0.61     1.39
		5        7      7.3           3            2.3           0.61     1.39
		6        7.5    0             0            0             0        0
	}\hybrid

	\pgfplotstableread{
		Location Force  GlobalSolve   InnerSolve   CoarseSolve   GMRES    Other
		1        4      9.1           3.38         3             0.96     1.76
		2        4.5    9.1           3.4          2.85          1        1.85
		3        5      9.2           3.5          3             1        1.7
		4        6      9.2           3.6          2.8           1        1.8
		5        7      0             0            0             0        0
		6        7.5    0             0            0             0        0
	}\RGDSWtwo
 

	\pgfplotstableread{
		Location Force  GlobalSolve   GMRES    Other
		1        4      5.1           4.4      0.7
		2        4.5    6.1           5.2      0.9
		3        5      0             0        0
		4        6      0             0        0
		5        7      0             0        0
		6        7.5    0             0        0
	}\NKS

	\begin{axis}[bar shift=0pt, hide axis]
		\node[mynodestyle]at (axis cs:1,0) {A};
		\node[mynodestyle] (two) at (axis cs:1,9.1) {$75$ $(4)$};
		\node[mynodestyle] at(axis cs:2,9.1) {$76$ $(4)$};
		\node[mynodestyle] at(axis cs:3,9.1) {$76$ $(4)$};
		\node[mynodestyle] at(axis cs:4,9.1) {$76$ $(4)$};
		\node at(axis cs:5,.4) {\scriptsize\color{red}\ding{55}};
		\node at(axis cs:6,.4) {\scriptsize\color{red}\ding{55}};

		\addplot+ table [x=Location, y=InnerSolve] {\RGDSWtwo};
		\addplot+ table [x=Location, y=CoarseSolve] {\RGDSWtwo};
		\addplot+ table [x=Location, y=GMRES] {\RGDSWtwo};
		\addplot+ table [x=Location, y=Other] {\RGDSWtwo};
	\end{axis}

	\begin{axis}[bar shift=11pt]
		\node[mynodestyle] at ([xshift=11pt]axis cs:1,0) {H};
		\node[xshift=11pt,mynodestyle] (three) at (axis cs:1,6.9) {$44$ $(3)$};
		\node[xshift=11pt,mynodestyle] at (axis cs:2,7.3) {$44$ $(3)$};
		\node[xshift=11pt,mynodestyle] at (axis cs:3,6.8) {$44$ $(3)$};
		\node[xshift=11pt,mynodestyle] at (axis cs:4,7.1) {$45$ $(3)$};
		\node[xshift=11pt,mynodestyle] at (axis cs:5,7.3) {$47$ $(3)$};
		\node[xshift=11pt,rotate=0] at  (axis cs:6,.4) {\scriptsize\color{red}\ding{55}};

		\addplot+ table [y=InnerSolve] {\hybrid}; \addlegendentry{Inner solve}
		\addplot+ table [y=CoarseSolve] {\hybrid}; \addlegendentry{Coarse solve}
		\addplot+ table [y=GMRES] {\hybrid}; \addlegendentry{GMRES}
		\addplot+ table [y=Other] {\hybrid}; \addlegendentry{Other}
	\end{axis}

	\begin{axis}[bar shift=-11pt, hide axis, cycle list shift=2]
		\node[mynodestyle]at ([xshift=-11pt]axis cs:1,0) {NKS};
		\node[xshift=-11pt,mynodestyle] (one) at (axis cs:1,5.2) {$124$ $(5)$};
		\node[xshift=-11pt,mynodestyle] at (axis cs:2,6.2) {$133$ $(6)$};
		\node[xshift=-11pt,rotate=0] at (axis cs:3,.4) {\scriptsize\color{red}\ding{55}};
		\node[xshift=-11pt,rotate=0] at (axis cs:4,.4) {\scriptsize\color{red}\ding{55}};
		\node[xshift=-11pt,rotate=0] at (axis cs:5,.4) {\scriptsize\color{red}\ding{55}};
		\node[xshift=-11pt,rotate=0] at (axis cs:6,.4) {\scriptsize\color{red}\ding{55}};
		\node[text width=1.8cm] (gmres) at (axis cs:1,13.8) {GMRES its. (Newton its.)};

		\addplot+ table [x=Location, y=GMRES] {\NKS};
		\addplot+ table [x=Location, y=Other] {\NKS};
	\end{axis}

	\draw [thin] (gmres) --  (one);
	\draw [thin] (gmres) --  (two);
	\draw [thin] (gmres) --  (three);

\end{tikzpicture}
	\caption{Comparison of the robustness of \ac{NKS} and the two-level nonlinear Schwarz variants.}
	\label{fig:elasticity-nks-vs-ldc}
\end{figure}

Figure \ref{fig:elasticity-strong} illustrates the strong scaling behavior of the nonlinear Schwarz implementation. For comparison, the strong scaling of the \ac{NKS} solver is also shown. In this test configuration, the \ac{NKS} solver achieves faster runtimes. However, both nonlinear Schwarz variants exhibit the same strong scaling behavior as the \ac{NKS} solver, validating our implementation.
To compare how the nonlinear Schwarz solver handles relative nonlinearity with the Newton solver, we incrementally increased the volume load and solved the problem using \num{576} subdomains. The results are shown in Figure \ref{fig:elasticity-nks-vs-ldc}. While the \ac{NKS} solver is approximately $30\%$ faster in this configuration, both nonlinear Schwarz variants demonstrate better robustness than the \ac{NKS} solver, with the hybrid variant being the most robust of the two.

\subsection{Coarse space observations for lid-driven cavity flow and nonlinear elasticity}
\label{sec:coarse-observations}
In this section, we discuss key observations regarding the coarse space. Our tests revealed that the coarse space plays a crucial role in the performance of the two-level nonlinear Schwarz method. In particular, small changes in the coarse space have a significantly larger impact on the performance of the nonlinear Schwarz solver compared to its linear counterpart. We begin with the modification we made to the \ac{RGDSW}/\ac{MsFEM} coarse spaces, then compare the performance of different coarse spaces in solving the \ac{LDC} problem, and finally visualize the differences between the \ac{RGDSW} and \ac{MsFEM} coarse spaces constructed for the \ac{LDC} problem.

\begin{figure}[h!]
    \centering
	\input{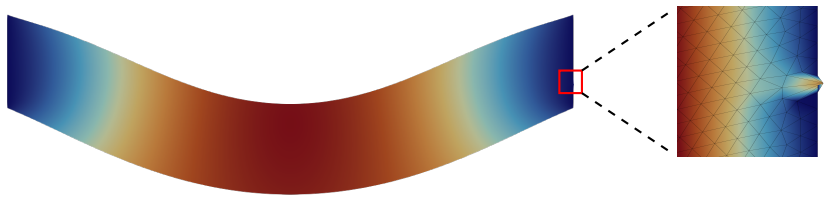}
	\caption{The beam model under an applied volume force of $9\;\unit{\MN/m^{2}}$.}
	\label{fig:beam-coarse-mesh}
\end{figure}

\begin{figure}[h!]
	\begin{subfigure}{0.32\textwidth}
		\centering
        \includegraphics[width=0.9\textwidth,height=2.5cm]{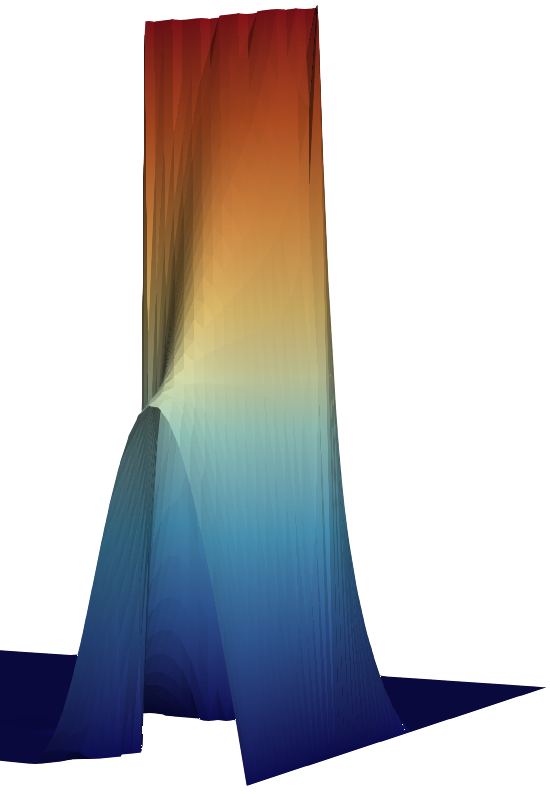}\\[3mm]
		\includegraphics[width=0.9\textwidth,height=2.5cm]{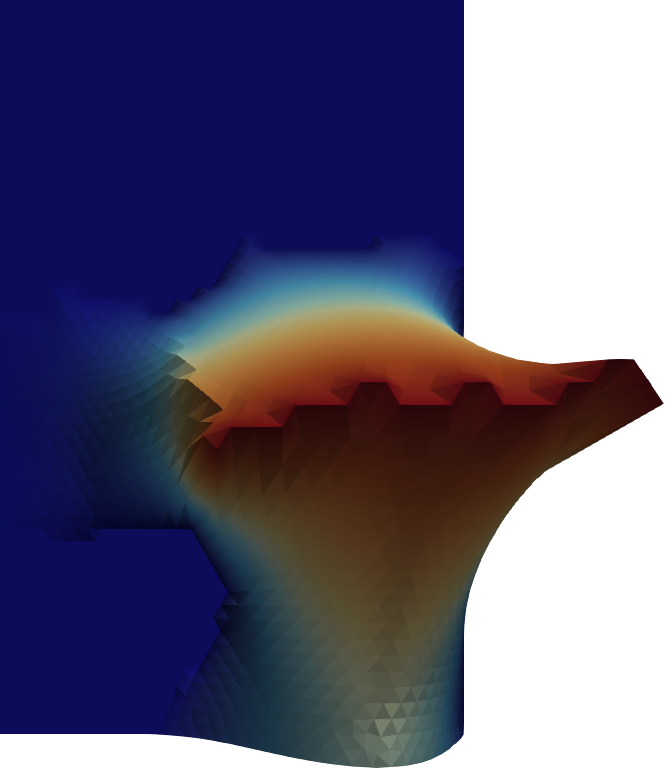}
		\caption{Unmodified.}
		\label{fig:beam-coarse-unmodified}
	\end{subfigure}
	\hfill%
	\begin{subfigure}{0.32\textwidth}
		\centering
        \includegraphics[width=0.9\textwidth,height=2.5cm]{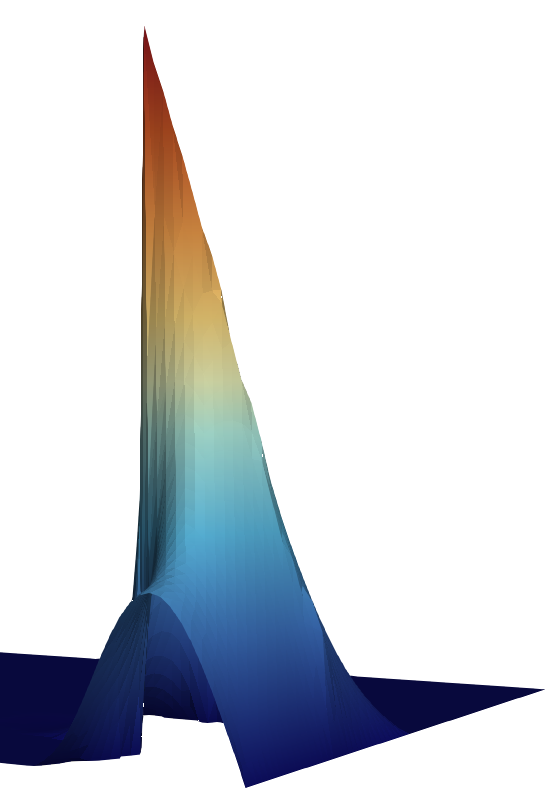}\\[3mm]
        \includegraphics[width=0.9\textwidth,height=2.5cm]{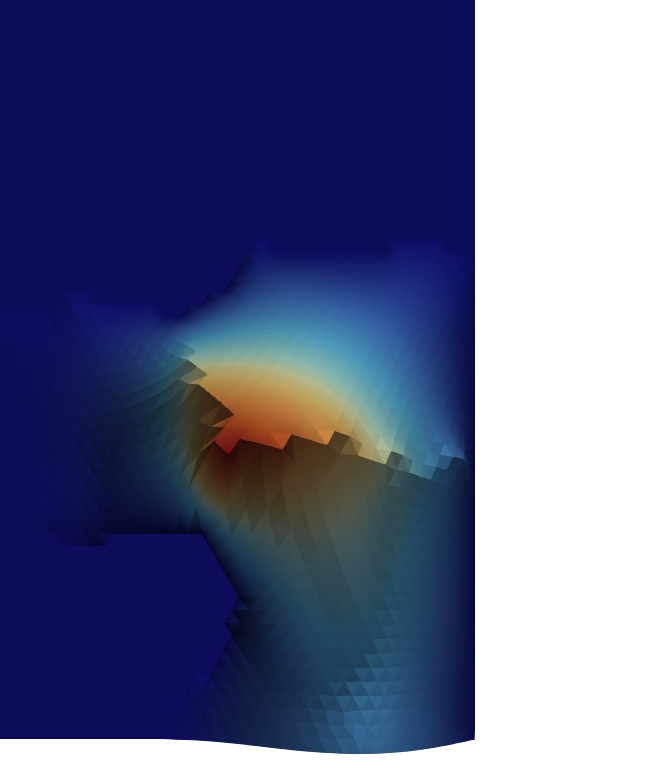}
		\caption{Modified.}
		\label{fig:beam-coarse-modified}
	\end{subfigure}
    \hfill%
	\begin{subfigure}{0.32\textwidth}
		\centering
        \includegraphics[width=0.6\textwidth,height=2.5cm]{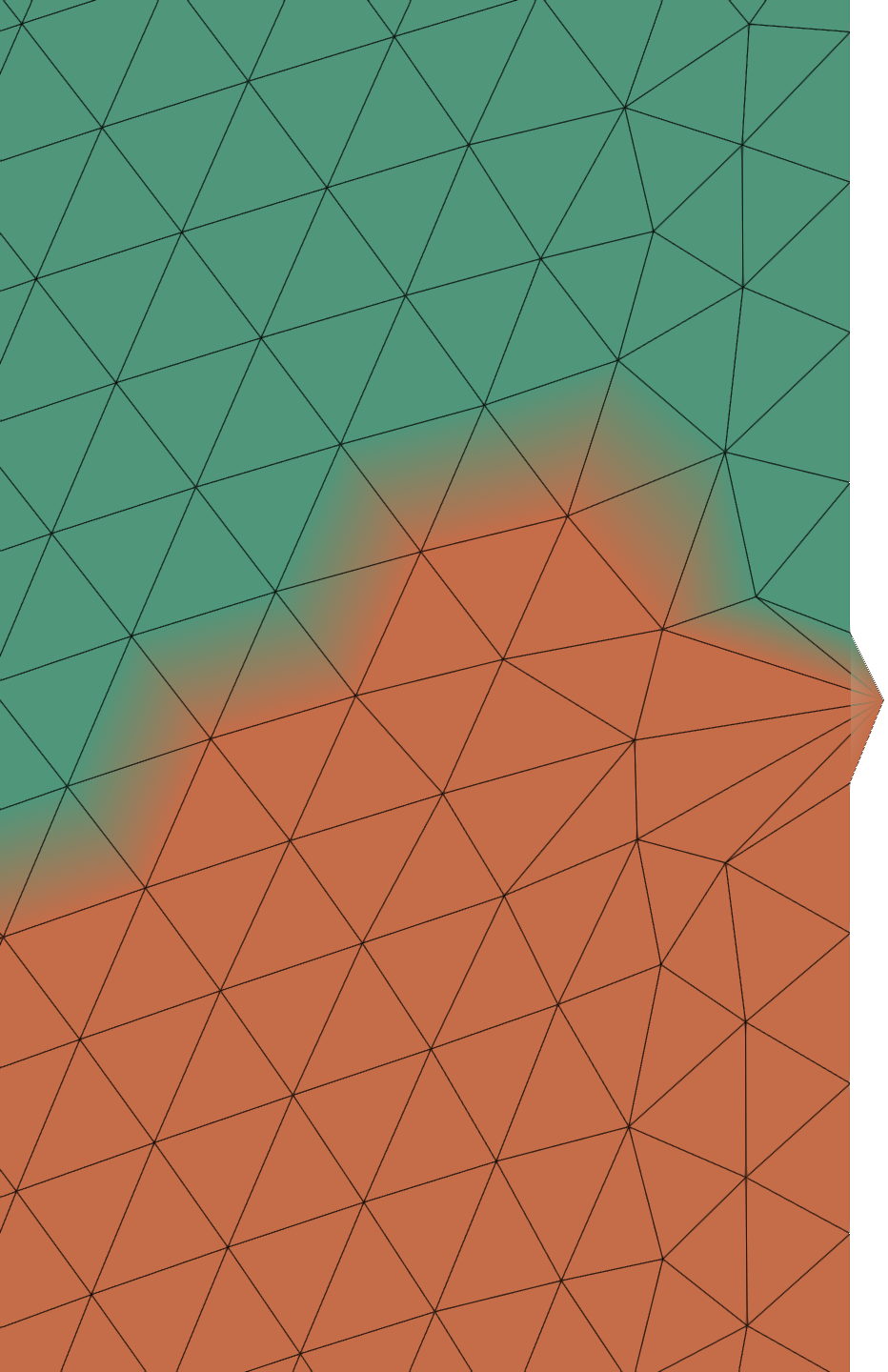}\\[3mm]
        \includegraphics[width=0.6\textwidth,height=2.5cm]{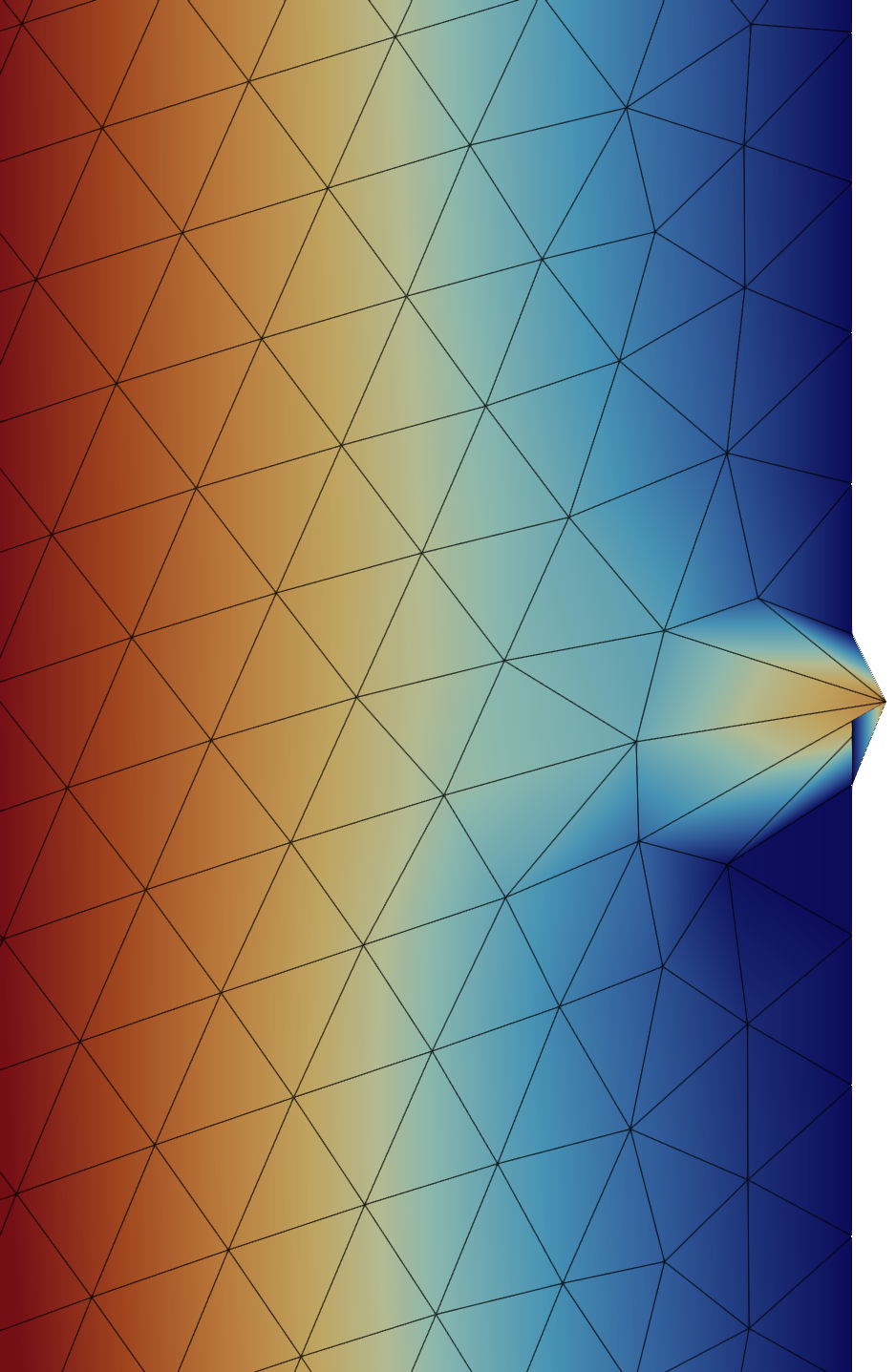}
		\caption{Colored for deformation.}
		\label{fig:beam-instability-deformation}
	\end{subfigure}
    \caption{The \ac{MsFEM} coarse basis for the nonlinear elasticity problem is constructed as described in Section~\ref{sec:nltwo}, with one coarse basis function per displacement direction at each interface vertex. Shown are different perspectives of the $y$-component of a coarse basis function corresponding to the $y$-direction at an interface edge ending at a Dirichlet boundary. The mesh consists of \num{26e3} nodes and \num{72} subdomains.}
	\label{fig:beam-coarse-basis}
\end{figure}

In our initial tests on the nonlinear elasticity problem \eqref{eq:nonlinelas}, the nonlinear Schwarz solver did not perform well using the \ac{MsFEM} coarse space described in Section \ref{sec:coarse}. Figure \ref{fig:beam-coarse-mesh} shows the deflected beam for a large volume force, this time using a mesh with \num{26e3} nodes. The nonlinear solvers are more easily able to solve the problem resulting from the coarser mesh, allowing a higher volume force to be used than in the tests discussed in Section \ref{sec:2dbeam}. A closer look at the right edge of the beam shows that a node of the mesh has crossed the zero Dirichlet boundary. The node lies directly on the interface edge between two subdomains. In Figure \ref{fig:beam-coarse-unmodified} we visualize the $y$-component of the $y$ degree of freedom of the \ac{MsFEM} coarse basis function corresponding to the interface node adjacent to the interface edge mentioned above. The coarse basis function has a constant value of one at the interface edge and drops to zero at the Dirichlet boundary within a single element. This is by design in the \ac{MsFEM} coarse space at interface edges ending at a Dirichlet boundary to maximize the area where the coarse space gives a partition of unity. The node shown in Figure \ref{fig:beam-instability-deformation} which has passed through the Dirichlet boundary, is exactly the last node on the interface where the coarse basis function has a value of one. This suggests that the deformation in the $y$ direction at the interface node, which dictates the deformation for the entire coarse basis function, is large enough to cause the observed unphysical effect.

To remedy this, we modified the construction of the \ac{MsFEM} coarse basis functions in \ac{FROSch} to have an inverse Euclidean decrease along the interface edges ending at a Dirichlet boundary, mirroring the way interface values are determined at internal interface edges as described by \eqref{eq:inveuclidean}. The modified coarse basis function is shown in Figure \ref{fig:beam-coarse-modified} on the same edge shown in the other figures. The modification solves the problem because the large deformation in the $y$ direction at the interface node is not propagated to all nodes in the interface edge. Modifying both the \ac{MsFEM} and \ac{RGDSW} coarse spaces in this way significantly improved the robustness of the nonlinear Schwarz solver to handle highly nonlinear problems, leading us to use it in our tests discussed in Sections \ref{sec:ldc} and \ref{sec:2dbeam}.

\begin{figure}[h!]
	\begin{tikzpicture}
	\pgfplotsset{
		every axis/.append style={
				ybar stacked,
				width=\textwidth,
				height=0.5\textwidth,
				ylabel={Runtime (seconds)},
				xlabel={$Re$},
				symbolic x coords={500, 750, 1000, 1500, 2000},
				xtick=data,
				enlarge x limits=0.1,
				legend style={at={(1,1)},anchor=north east},
				axis lines*=left, ymajorgrids, yminorgrids,
				ymin=0,
				ymax=200,
				bar width=8pt,
				minor y tick num=1,
				xticklabel style={rotate=0,xshift=0ex,anchor=north},
				cycle list name=Set2-5,
			},
		every axis plot/.append style={
				fill,
			},
	}
	\begin{axis}[bar shift=-11pt, hide axis]
		\node [rotate=90](rgdsw) at ([xshift=-11pt]axis cs:500,25) {RGDSW};
		\addplot coordinates {(500,22) (750,23.7) (1000,23.8) (1500,31) (2000,0)};
		\addplot coordinates {(500,17.4) (750,17.1) (1000,18.3) (1500,24) (2000,0)};
		\addplot coordinates {(500,28.6) (750,37) (1000,42.2) (1500,64) (2000,0)};
		\addplot coordinates {(500,23.7) (750,22.6) (1000,23.7) (1500,28) (2000,0)};

		\node [rotate=90,anchor=center](500) at ([xshift=-11pt]axis cs:500,115) {$320$ $(4)$};
		\node [rotate=90,anchor=center](750) at ([xshift=-11pt]axis cs:750,123) {$400$ $(4)$};
		\node [rotate=90,anchor=center](1000) at ([xshift=-11pt]axis cs:1000,132) {$450$ $(4)$};
		\node [rotate=90,anchor=center](1500) at ([xshift=-11pt]axis cs:1500,170) {$650$ $(5)$};
		\node at ([xshift=-11pt]axis cs:2000,6) {\scriptsize\color{red}\ding{55}};

	\end{axis}

	\begin{axis}
		\node [rotate=90](msfem) at (axis cs:500,24) {MsFEM};
		\addplot coordinates {(500,21.9) (750,21.9) (1000,0) (1500,0) (2000,0)};
		\addplot coordinates {(500,22.3) (750,24  ) (1000,0) (1500,0) (2000,0)};
		\addplot coordinates {(500,25.6) (750,28.9) (1000,0) (1500,0) (2000,0)};
		\addplot coordinates {(500,23.7) (750,23.9) (1000,0) (1500,0) (2000,0)};

		\legend{
			Inner solve,
			Coarse solve,
			GMRES,
			Other
		}
		\node[rotate=90,anchor=center](one) at (axis cs:500,117) {$300$ $(4)$};
		\node[rotate=90,anchor=center] at (axis cs:750,122) {$330$ $(4)$};

		\node at (axis cs:1000,6) {\scriptsize\color{red}\ding{55}};
		\node at (axis cs:1500,6) {\scriptsize\color{red}\ding{55}};
		\node at (axis cs:2000,6) {\scriptsize\color{red}\ding{55}};

	\end{axis}
    
	\begin{axis}[bar shift=11pt, hide axis]
		\node [rotate=90](gdsw) at ([xshift=11pt]axis cs:500,22) {GDSW};
		\addplot coordinates {(500,25.3) (750,0) (1000,0) (1500,0) (2000,0)}; 
		\addplot coordinates {(500,23.7) (750,0  ) (1000,0) (1500,0) (2000,0)}; 
		\addplot coordinates {(500,32.4) (750,0) (1000,0) (1500,0) (2000,0)}; 
		\addplot coordinates {(500,29.2) (750,0) (1000,0) (1500,0) (2000,0)}; 

		\node [rotate=90,anchor=center](two) at ([xshift=11pt]axis cs:500,133) {$370$ $(5)$};
		\node at ([xshift=11pt]axis cs:750,6) {\scriptsize\color{red}\ding{55}};
		\node at ([xshift=11pt]axis cs:1000,6) {\scriptsize\color{red}\ding{55}};
		\node at ([xshift=11pt]axis cs:1500,6) {\scriptsize\color{red}\ding{55}};
		\node at ([xshift=11pt]axis cs:2000,6) {\scriptsize\color{red}\ding{55}};

	\end{axis}

	\node[rotate=0, text width=1.8cm] (gmres) at ([xshift=11,yshift=50]500){GMRES its. (Newton its.)};

	\draw [thin] (gmres) --  (500);
	\draw [thin] (gmres) --  (one);
	\draw [thin] (gmres) --  (two);


\end{tikzpicture}
	\caption{Time-to-solution of the hybrid solver for the LDC problem using the GDSW and the modified RGDSW and MsFEM coarse spaces. \num{256} subdomains were used and the tolerances were set as in Section \ref{sec:ldc}.}
	\label{fig:rgdsw-vs-msfem-vs-gdsw-ldc}
\end{figure}
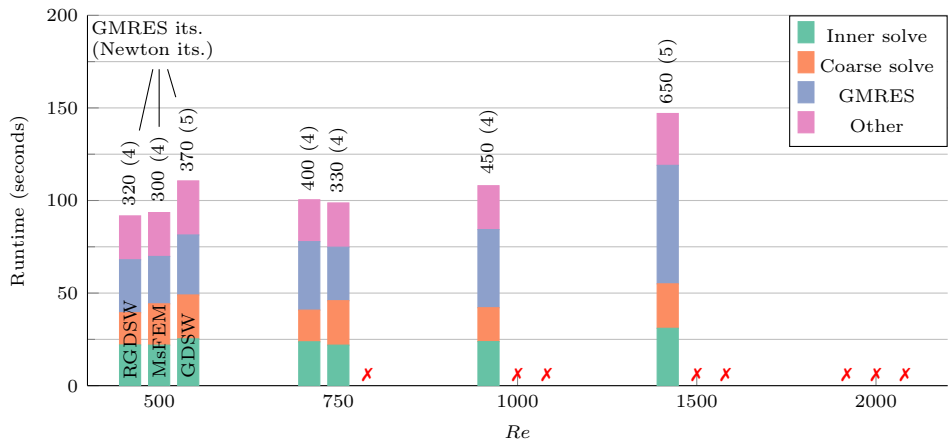

The modified coarse space produced satisfactory results for the nonlinear elasticity problem \eqref{eq:nonlinelas}, which served as a proof of concept. Thus, we did not explore alternative coarse spaces for this case. In contrast, the application of our nonlinear Schwarz implementation to the \ac{LDC} problem required some numerical testing, as the selection of an appropriate coarse space was more challenging and a guiding theory does not exist. Table \ref{tab:coarse-comparison} compares the convergence performance of different coarse spaces. The modified variants generally outperform their unmodified counterparts, with the modified \ac{MsFEM} appearing to be the most effective for this problem. Since the \ac{GDSW} coarse space has individual coarse basis functions on interface edges ending at Dirichlet nodes, we would expect it to perform in the same ballpark as the modified versions of the other coarse spaces. This is confirmed by the results in Table \ref{tab:coarse-comparison}. Figure \ref{fig:rgdsw-vs-msfem-vs-gdsw-ldc} shows time-to-solution results for the two-level nonlinear Schwarz solver using the \ac{GDSW} and modified coarse spaces over various Reynolds numbers and for larger subdomains. In this setting, the \ac{RGDSW} coarse space results in the most robust two-level nonlinear Schwarz method. Using the \ac{MsFEM} coarse space results in fewer \ac{GMRES} iterations, but the nonlinear solver fails to converge above $Re=750$. The longer runtime when using \ac{MsFEM}, despite the lower number of \ac{GMRES} iterations, is caused by an increase in the number of backtracking steps taken by the coarse Newton method, negating the \ac{GMRES} solver improvements. Since the \ac{GDSW} coarse space results in a slower, less robust method than both other coarse spaces, we used the modified \ac{RGDSW} coarse space for the tests in Section \ref{sec:ldc}. Recall that the robustness of our solver increases with the number of subdomains, such that we were able to solve for higher Reynolds numbers than shown in Figure \ref{fig:rgdsw-vs-msfem-vs-gdsw-ldc}. We suppose that this property is preserved when using other coarse spaces, however, we have not tested this.

The idea of mixing coarse space types in monolithic two-level linear Schwarz preconditioners is explored in \cite{heinleinMonolithicBlockOverlapping2025, sassmanshausenDiss}. Specifically for fluid problems, a combination of \ac{GDSW} for the velocity degrees of freedom and \ac{RGDSW} for the pressure degrees of freedom performs well. With some modifications, this idea could be applied using our nonlinear Schwarz implementation. We reserve this for future research.

\begin{table}[!ht]
	\centering
	\renewcommand{\arraystretch}{1.2} 
	{\small \begin{tabular}{|p{2.3cm}|>{\raggedleft\arraybackslash}p{1.7cm}|>{\raggedleft\arraybackslash}p{1.7cm}|>{\raggedleft\arraybackslash}p{0.9cm}|>{\raggedleft\arraybackslash}p{0.9cm}|>{\raggedleft\arraybackslash}p{1.1cm}|>{\raggedleft\arraybackslash}p{1.3cm}|}
			\hline
			\multicolumn{1}{|l|}{}         &
			\multicolumn{1}{l|}{Rel. Res.} &
			\multicolumn{1}{l|}{Abs. Res.} &
			\multicolumn{1}{l|}{Outer}     &
			\multicolumn{1}{l|}{Inner}     &
			\multicolumn{1}{l|}{Coarse}    &
			\multicolumn{1}{l|}{GMRES}                                                                        \\
			\hline
			Mod. RGDSW                     & $1 \times 10^{-6}$   & $1.6 \times 10^{-8}$ & 4 & 6   & 10 & 94  \\ \hline
			RGDSW                          & $1.6 \times 10^{-5}$ & $2.5 \times 10^{-7}$ & 4 & 6.7 & 12 & 100 \\ \hline
			Mod. MsFEM                     & $3.8 \times 10^{-7}$ & $5.8 \times 10^{-9}$ & 4 & 5   & 16 & 79  \\ \hline
			MsFEM                          & -                    & -                    & - & -   & -  & -   \\ \hline
			GDSW                           & $2.9 \times 10^{-5}$ & $4.4 \times 10^{-7}$ & 4 & 7.1 & 14 & 81  \\ \hline
		\end{tabular}}
	\caption{Performance comparison of various coarse spaces with the hybrid variant of the two-level nonlinear Schwarz solver for the LDC problem at $Re=1000$ with \num{256} subdomains. Small subdomains with $H/h=15$ were used, an overlap of 5 was set, and solver tolerances are chosen as described in Section \ref{sec:ldc}. Modified coarse spaces are labeled \textit{mod.} The columns \textit{Outer}, \textit{Inner}, and \textit{Coarse} show the respective Newton iterations, with inner and GMRES iterations summed over all outer iterations and inner iterations averaged across subdomains. Missing entries indicate that the method failed to converge.}
	\label{tab:coarse-comparison}
\end{table}

\section*{Acknowledgments}
The authors gratefully acknowledge the financial support {\KH from} the German Federal Ministry of Research, Technology and Space (BMFTR) in the program SCALEXA.
The scientific support and HPC resources provided by the Erlangen National High Performance Computing Center (NHR@FAU) of the Friedrich-Alexander-Universit\"at Erlangen-N\"urnberg (FAU), Germany under the NHR project k107ce is gratefully acknowledged. NHR funding is provided by federal and Bavarian state authorities. NHR@FAU hardware is partially funded by the German Research Foundation (DFG) - 440719683. The authors would like to thank Alexander Heinlein (TU Delft) for his help with \ac{FROSch} implementation details, Lea Saßmannshausen (University of Cologne) for her help using the \ac{FEDDLib} software library, and Sharan Nurani Ramesh (Ruhr University Bochum) for his help with the hyperelasticity problem. During the preparation of this work the generative AI tools ChatGPT, Claude, and Apple Intelligence were used by the first author to correct grammar and spelling and to improve style of writing. After using these tools, all authors reviewed and edited the content as needed and take full responsibility for the content of the publication.

\bibliographystyle{siamplain} 
\bibliography{PhD}

\end{document}